\documentclass{amsart}

\usepackage{amscd,amssymb}

\usepackage[all]{xy}

\renewcommand{\theenumi}{\roman{enumi}}

\def\Z{{\mathbb Z}}
\def\Q{{\mathbb Q}}
\def\R{{\mathbb R}}
\def\C{{\mathbb C}}

\def\F{{\mathbb F}}
\def\P{{\mathbb P}}
\def\L{{\mathbb L}}
\def\H{{\mathbb H}}
\def\V{{\mathbb V}}

\def\J{{\mathcal J}}
\def\K{{\mathcal K}}
\def\M{{\mathcal M}}

\def\U{{\mathcal U}}
\def\cC{{\mathcal C}}
\def\cE{{\mathcal E}}
\def\cG{{\mathcal G}}
\def\cH{{\mathcal H}}

\def\g{{\mathfrak g}}
\def\p{{\mathfrak p}}
\def\u{{\mathfrak u}}

\def\l{{\ell}}

\def\e{{\epsilon}}
\def\w{{\omega}}
\def\k{{\kappa}}
\def\G{{\Gamma}}

\def\Bbar{\overline{B}}
\def\Cbar{\overline{C}}

\def\Ghat{\widehat{\G}}
\def\Kbar{\overline{K}}

\def\Wtilde{\widetilde{W}}
\def\Xbar{\overline{X}}

\def\rhotilde{{\tilde{\rho}}}
\def\tautilde{{\tilde{\tau}}}
\def\etilde{{\tilde{\epsilon}}}
\def\rtilde{{\tilde{r}}}
\def\rhohat{{\hat{\rho}}}
\def\ctilde{{\tilde{c}}}
\def\mhat{{\hat{m}}}
\def\nhat{{\hat{n}}}
\def\utilde{{\tilde{u}}}
\def\xihat{{\hat{\xi}}}
\def\xtilde{{\tilde{x}}}
\def\nutilde{{\tilde{\nu}}}

\def\un{\mathrm{un}}
\def\et{{\mathrm{\acute{e}t}}}
\def\alg{{\mathrm{alg}}}
\def\hodge{{\mathrm{Hodge}}}
\def\hom{{\mathrm{hom}}}

\def\prol{{(\ell)}}

\def\an{\mathrm{an}}

\def\Ql{{\Q_\ell}}
\def\Zl{{\Z_\ell}}
\def\Zlx{{\Z_\ell^\times}}

\def\Het{H_\et}		
\def\Gm{{\mathbb{G}_m}}
\def\GSp{{\mathrm{GSp}}}
\def\Sp{{\mathrm{Sp}}}

\def\Cyc{{\mathrm{Cyc}}}

\def\Can{{C^\an}}
\def\Canx{{C^\an,\xi^\an}}
\def\Cx{{C,\xi}}
\def\xian{{\xi^\an}}

\def\dot{{\bullet}}

\def\blank{\phantom{x}}
\def\comp{~\widehat{\!}{\;}}

\newcommand\im{\operatorname{im}}
\newcommand\id{\operatorname{id}}
\newcommand\Spec{\operatorname{Spec}}
\newcommand\Hom{\operatorname{Hom}}
\newcommand\End{\operatorname{End}}
\newcommand\Ext{\operatorname{Ext}}
\newcommand\Aut{\operatorname{Aut}}

\newcommand\Out{\operatorname{Out}}
\newcommand\Inn{\operatorname{Inn}}
\newcommand\Gal{\operatorname{Gal}}
\newcommand\Gr{\operatorname{Gr}}
\newcommand\Diff{\operatorname{Diff}}
\newcommand\Jac{\operatorname{Jac}}
\newcommand\Lib{\operatorname{Lib}}

\newtheorem{theorem}{Theorem}[section]
\newtheorem{lemma}[theorem]{Lemma}
\newtheorem{proposition}[theorem]{Proposition}
\newtheorem{corollary}[theorem]{Corollary}
\newtheorem{bigtheorem}{Theorem}

\theoremstyle{definition}
\newtheorem{definition}[theorem]{Definition}

\theoremstyle{remark}
\newtheorem{remark}[theorem]{Remark}

\begin{document}

\title[Galois actions on fundamental groups and $C-C^-$]
{Galois actions on fundamental groups of curves and the cycle $C-C^-$}

\author[Richard Hain]{Richard Hain}
\address{Department of Mathematics\\ Duke University\\
Durham, NC 27708-0320}
\email{hain@math.duke.edu}
\thanks{Supported in part by grant DMS-0103667 from the National Science
Foundation.}

\author[Makoto Matsumoto]{Makoto Matsumoto}
\address{Department of Mathematics \\ Hiroshima University\\
Hiroshima, 739-8526 JAPAN}
\email{m-mat@math.sci.hiroshima-u.ac.jp}

\date{\today}



\maketitle

\section{Introduction}
\label{sec:introduction}

The goal of this paper is to better understand the action of Galois groups on
fundamental groups of smooth curves of genus $3$ or more, and to show, using a
fundamental result of Dennis Johnson \cite{johnson:h1}, how the algebraic cycle
$C-C^-$ in the jacobian of a curve $C$ helps control the size of the image of
the Galois group in the automorphism group of the curve's fundamental group.

Suppose that $K$ is a field of characteristic zero and that $X$ is a smooth
projective variety defined over $K$. Denote $\Gal(\Kbar/K)$ by $G_K$. Fix a
prime number $\l$. By standard constructions (Section~\ref{sec:Abel-Jacobi}), a
codimension $r$ algebraic cycle $Z$ in $X$, defined over $K$ and homologically
trivial in $X\otimes \Kbar$, determines a class
$$
e_Z \in H^1(G_K,\Het^{2r-1}(X\otimes \Kbar,\Zl(r))).
$$
This depends only on the rational equivalence class of $Z$ and is the image
of $Z$ under the $\l$-adic Abel-Jacobi mapping
$$
CH^r_\hom(X) \to H^1(G_K, \Het^{2r-1}(X \otimes \Kbar, \Zl(r)))
$$
which is defined on the group of rational equivalence classes of
homologically trivial, codimension $r$ cycles on $X$.

Suppose that $C$ is a smooth, geometrically connected projective curve over $K$
of genus $g \ge 3$ and that $\xi$ is a $K$-rational point of $C$. The morphism
$\sigma_\xi : C \to \Jac C$ that takes $z \in C$ to the divisor class of
$z-\xi$ is an embedding. Define the algebraic 1-cycle $C_\xi$ in $\Jac C$ to be
$(\sigma_\xi)_\ast C$. One also has the cycle $C_\xi^- := i_\ast C_\xi$, where
$i$ is the involution of the jacobian that takes each point to its inverse.

Two algebraic cycles particularly relevant to Galois actions on the fundamental
group of the algebraic curve $(C,\xi)$ are:
\begin{enumerate}
\item the $0$-cycle $(2g-2)\xi - K_C$ in $C$, where $K_C$ is any canonical
divisor of $C$;
\item the 1-cycle $C_\xi - C_\xi^-$ in $\Jac C$.
\end{enumerate}
Both are homologically trivial. The first cycle defines a class
$$
\k(\Cx) \in H^1(G_K,\Het^1(C \otimes \Kbar ,\Zl(1)))
$$
and the second a class
$$
\mu(\Cx) \in H^1(G_K,\Het^{2g-3}(\Jac C \otimes \Kbar,\Zl(g-1))).
$$

It is convenient to set $H_\Zl = \Het^1(C\otimes \Kbar,\Zl(1))$  and $L_\Zl =
(\Lambda^3 H_\Zl)(-1)$. Both are of weight $-1$.  Denote their tensor products
with $\Ql$ by $H_\Ql$ and $L_\Ql$,  respectively. Wedging with the polarization
$q \in \Lambda^2 H_\Zl (-1)$  defines a $G_K$-invariant embedding $H_\Zl \to
L_\Zl$. Define
$$
\nu(C) = \text{ the image of $\mu(\Cx)$ in } H^1(G_K,L_\Zl/H_\Zl).
$$
This is independent of the choice of $\xi$ and is defined, even when $C$ has no
$K$-rational points (cf.\ Section~\ref{sec:non-pointed}).

Suppose now that $K$ is a subfield of $\C$.\footnote{It is not really necessary
to assume that $K$ is a subfield of $\C$, but we have chosen to do this in the
hope of making the paper a little more accessible to non-experts. Experts
should have no trouble generalizing our main results to arbitrary fields $K$ of
characteristic zero and using standard arguments to deduce them from statements
proved in this paper.} Let $\xi^\an$ be the geometric point $\xi\otimes \C$ on
the complex analytic space $C^\an$ associated to $C$. Denote the Lie algebra of
the $\Ql$-form of the unipotent completion of the topological fundamental group
$\pi_1(\Canx)$ by $\p(\Cx)$. This is a pronilpotent Lie algebra over $\Ql$.

Our goal in this paper is to describe how the classes $\nu(C)$, $\mu(C,\xi)$
and $\kappa(\Cx)$ influence the size of the image of $G_K$ in the automorphism
group of $\p(\Cx)$. In particular, we are interested in the relationship
between the image of $G_K$ in $\Aut \p(\Cx)$ and the image of the mapping class
group.

We shall denote the mapping class group of a closed, pointed, oriented surface
$(S,x)$ by $\G_{S,x}$. It is the group of connected components of the group of
orientation preserving diffeomorphisms of $S$ that fix $x$:
$$
\G_{S,x} = \pi_0 \Diff^+ (S,x).
$$
There is a natural homomorphism $\G_{S,x} \to \Aut \pi_1(S,x)$. In the
un-pointed case we shall denote the mapping class group of $S$ by $\G_S$. It is
the group of connected components of the group of orientation preserving
diffeomorphisms of $S$:
$$
\G_S = \pi_0 \Diff^+ S.
$$
There is a natural homomorphism $\G_S \to \Out \pi_1(S,x)$.

Denote the lower central series of $\p(C,\xi)$ by
$$
\p(C,\xi) = L^1\p(\Cx) \supset L^2 \p(\Cx) \supset L^3 \p(\Cx) \supset
\cdots .
$$
The group $\Aut \p(\Cx)$ is an affine proalgebraic group, as it is the inverse
limit of the automorphism groups of the $\p(\Cx)/L^m$. The group
$\pi_1(C^\an,\xi^\an)$ acts on itself by conjugation, and, by functoriality, on
$\p(\Cx)$. Denote the Zariski closure of its image by $\Inn \p(\Cx)$. Set
$$
\Out \p(\Cx) = \Aut \p(\Cx)/\Inn \p(\Cx).
$$
This is also a proalgebraic group.

By functoriality, there are natural homomorphisms
$$
\theta_\Cx : \G_\Canx \to \Aut \p(\Cx)\text{ and }
\theta_C : \G_\Can \to \Out \p(\Cx).
$$
Since the filtration $L^\dot$ is by characteristic subalgebras, we also have
the $m$-truncated representations
$$
\theta_\Cx^m : \G_\Canx \to \Aut \big[\p(\Cx)/L^{m+1}\big] \text{ and }
\theta_C^m : \G_\Can \to \Out \big[\p(\Cx)/L^{m+1}\big].
$$
When $m=1$, these are the standard actions of $\G_\Cx$ and $\G_C$ on the first
homology of $\Can$.

On the other hand, the theory of algebraic fundamental groups gives nonabelian
Galois representations
$$
\rho^\prol_\Cx : G_K \to \Aut \pi_1^\prol(\Canx)
\text{ and }
\rho^\prol_C : G_K \to \Out \pi_1^\prol(\Canx),
$$
where $\pi_1^\prol(\Canx)$ denotes the pro-$\l$ completion of $\pi_1(\Canx)$.
These induce homomorphisms
$$
\rho_\Cx : G_K \to \Aut \p(\Cx) \text{ and } \rho_C : G_K \to \Out \p(\Cx).
$$
We also have their $m$-truncated versions
$$
\rho_\Cx^m : G_K \to \Aut \big[\p(\Cx)/L^{m+1}\big] \text{ and }
\rho_C^m : G_K \to \Out \big[\p(\Cx)/L^{m+1}\big].
$$
When $m=1$, these two homomorphisms both coincide with the Galois action on the
abelianization of $\p(\Cx)$, which is canonically isomorphic to $H_\Ql$. Since
this action preserves the cup product up to a character, the image of
$\rho_C^1$ is contained in the subgroup $\GSp(H_\Ql)$ of $\Aut H_\Ql$. We will
view $\rho_C^1$ as a homomorphism
$$
\rho_C^1: G_K \to \GSp(H_\Ql).
$$

The following theorem explains how the representation $\rho_C^1$ and the
classes $\nu(C)$ and $\kappa(\Cx)$ determine when the Zariski closure of the
image of $G_K$ in $\Aut \p(\Cx)$ contains the geometric automorphisms $\im
\theta_\Cx$. Recall that $C$ has genus $\ge 3$.

\begin{bigtheorem}
\label{th:main-aut}
If the $\l$-adic cyclotomic character $G_K \to \Zlx$ has infinite image, then
the following three conditions are equivalent:
\begin{enumerate}
\item The Zariski closure of the image of $\rho_\Cx$ contains the image of
$\theta_\Cx$.
\item  The Zariski closure of the image of $\rho^2_\Cx$ in
$\Aut\big[\p(\Cx)/L^3\big]$ contains the image of $\theta_\Cx^2$.
\item The homomorphism $\rho^1_C : G_K \to \GSp(H_\Ql)$ is Zariski dense, and
both the classes $\k(\Cx)\in H^1(G_K,H_\Zl)$  and $\nu(C)\in
H^1(G_K,L_\Zl/H_\Zl)$ have infinite order.
\end{enumerate}
\end{bigtheorem}

\begin{remark} Several remarks are in order:
\begin{enumerate}
\item[(a)] This theorem has an alternative formulation,
Theorem~\ref{th:main-aut-4}, in terms of pro-$\l$ completions instead of
unipotent completions.

\item[(b)] If one (and consequently all) of the conditions in the theorem are
satisfied, then the Zariski closure  of the image of $G_K$ in $\Aut \p(\Cx)$
depends only on $g$, i.e., is independent of the choice of $(\Cx)$.

\item[(c)] Matsumoto and Tamagawa \cite[Thm.~1.2]{MAPGAL}  prove a result which
implies that there is an infinite number of genus $g$ curves defined over a
number field that satisfy condition (i). Thus both classes $\k(\Cx)$, $\nu(C)$
are of infinite order for an infinite number of genus $g$ curves, a fact that
can also be proved using the methods of Bloch-Esnault \cite{BE}. Note that it
is not known whether there is a number field $K$ for which there is an infinite
number of curves with this property.

\item[(d)] One can use any homologically trivial cycle defined for each smooth
curve $C$ whose normal function is a non-zero multiple of that of $C-C^-$. For
example, one can use the Gross-Schoen cycle \cite{gross-schoen}, which is a
homologically trivial 1-cycle in the triple product of a curve of genus $\ge 3$
whose class in $H^1(G_K,L_\Zl)$ is $3\nu(C)$.

\item[(e)] The theorem as stated does not hold in genus 2. However, there is a
similar theorem for hyperelliptic curves (where $\xi$ is a Weierstrass point)
which will be the subject of a future paper by the first author. The size of
the Galois image in the hyperelliptic case is controlled by two cohomology
classes which correspond to elements of Bloch's higher Chow groups. One part of
the analogue of the Harris-Pulte Theorem needed in this case is the result of
Colombo \cite{colombo}.
\end{enumerate}

\end{remark}

There is a similar result in the non-pointed case. In this case, the size of
the Zariski closure of the image of $\rho_C$ is controlled by $\nu(C)$ and
$\rho_C^1$. Note that the outer representations $\theta_C^m$ can be defined
even when $C$ has no $K$-rational points (cf.\ Section~\ref{sec:non-pointed}).
So for this theorem, we need only assume that $C$ is a smooth projective curve
of genus $\ge 3$ defined over the subfield $K$ of $\C$.

\begin{bigtheorem}
\label{th:main-out}
If the $\l$-adic cyclotomic  character $G_K \to \Zlx$ has infinite image, then
the following three conditions are equivalent:
\begin{enumerate}
\item The Zariski closure of the image of $\rho_C$ contains the
image of $\theta_C$.
\item The Zariski closure of the image of $\rho_C^2$ in
$\Out\big[\p(\Cx)/L^3\big]$ contains the image of $\theta_C^2$.
\item The homomorphism $\rho_C^1: G_K \to \GSp(H_\Ql)$ is Zariski dense and the
class $\nu(C)$ is of infinite order in $H^1(G_K,L_\Zl/H_\Zl)$.
\end{enumerate}
\end{bigtheorem}

Closely related to the density theorems above is the following $\l$-adic
analogue of the Harris-Pulte Theorem \cite[p.~722]{pulte}, which we prove
in Section~\ref{sec:harris-pulte}.

In Section~\ref{sec:char_class}, we show that the continuous action of a
topological group $\G$ on the pro-$\l$ fundamental group $\pi_1^\prol(\Canx)$
of $(\Canx)$ determines a characteristic class in the continuous cohomology
group $H^1(\G, L_\Zl)$. In particular, the action of $G_K$ on the pro-$\l$
fundamental group of $(\Canx)$ gives rise to a class
$$
m(\Cx) \in H^1(G_K, L_\Zl).
$$

\begin{bigtheorem}
\label{th:l-Harris-Pulte}
If the genus of $C$ is $\ge 3$, then
$$
\mu(\Cx) = m(\Cx) \in H^1(G_K, L_\Zl).
$$
\end{bigtheorem}

The original Harris-Pulte Theorem relates the mixed Hodge structure on
$\pi_1(\Cx)$ modulo the third term $L^3$ of its lower central series, where $C$
is a compact Riemann surface of genus $\ge 3$, to the element $\zeta(\Cx)$ of
$$
\Ext^1_\hodge(\Z,\Lambda^3 H^1(C^\an,\Z)(1))
$$
determined by $C_\xi - C_\xi^-$. The formulation \cite[Thm.~11.1]{hain:normal}
asserts that $\zeta(\Cx)=2\rho(\Cx)$, where $\rho(\Cx)$ is the extension class
of the mixed Hodge structure on $\pi_1(\Canx)/L^3$. Note that $2\rho(\Cx)$
corresponds to $m(\Cx)$.

An un-pointed version of Theorem~\ref{th:l-Harris-Pulte} is stated in
Section~\ref{sec:harris-pulte_unptd}.

We shall prove these three theorems by specialization from the case of the
universal curve over the moduli space $\M_g^1$ of pointed smooth projective
curves of genus $g$. A basic tool is the weighted completion of topological
groups, which we review in Section~\ref{sec:completion}.

\section{Notation, Conventions and Preliminaries}

Throughout $K$ will be a field of characteristic zero. In this article, all
schemes, varieties and stacks will be defined (and of finite type) over $K$
unless otherwise noted. By a stack, we shall mean a Deligne-Mumford stack
\cite{DM}. The category of schemes over a base scheme $B$ is a full subcategory
of the category of stacks over $B$.

Often $K$ will be a subfield of $\C$. In this case, we shall denote the
algebraic closure of $K$ in $\C$ by $\Kbar$ and the analytic variety
corresponding to a variety (resp.\ stack) $X$ over $K$ by $X^\an$. It will be
viewed as a topological space (resp.\ orbifold) in the complex topology. Each
complex point $\xi$ of $X$ determines a point of $X^\an$, which we shall
denote by $\xi^\an$.

The algebraic fundamental group of a stack (or scheme) $X$ with base point the
geometric point $x$ will be denoted by $\pi_1^\alg(X,x)$. We shall use the
notation $\pi_1(M,x)$ to denote the fundamental group of the pointed
topological space (or orbifold) $(M,x)$.

The profinite completion of a discrete group $\G$ will be denoted by $\Ghat$ or
$\G\comp$. The pro-$\l$ completion of a discrete or profinite group $\G$ will
be denoted by $\G^\prol$. The pro-$\l$ fundamental group functor will be
denoted by $\pi_1^\prol$ in both the algebraic and topological cases.

Suppose that $F$ is a field of characteristic zero. The unipotent (or Malcev)
completion over $F$ of the discrete group $\G$ will be denoted by
$\G^\un_{/F}$. If no field is mentioned, we will take $F$ to be $\Q$. The
prounipotent fundamental group functor will be denoted by $\pi_1^\un$.

Now suppose that $\l$ is a rational prime and that $F$ is an $\l$-adic field.
The continuous unipotent completion over $F$ of the profinite group $\G$
(\cite[Appendix~A.2]{hain-matsumoto}) will be denoted by $\G^\un_{/F}$. If no
field is mentioned, we will take $F$ to be $\Ql$. In
\cite[Thm.~A.6]{hain-matsumoto} it is shown that if $\G$ is a finitely
generated discrete group, then there are natural isomorphisms
$$
\G^\un_{/\Q} \otimes \Ql \cong  \G^\un_{/\Q_l}\cong 
\big(\G^\prol\big)^{\un}_{/\Ql} \cong \big(\Ghat\big)^{\un}_{/\Ql}
$$
of $\Ql$-prounipotent groups.

If $K$ is a subfield of $\C$, then for each variety $X$ over $K$ with an
$\C$-rational point $x$, there is, by the standard comparison theorems in
\cite{SGA1} and their generalizations, a canonical isomorphism
$$
\pi_1^\alg(X\otimes \Kbar,x) \cong \pi_1(X^\an,x^\an)\comp.
$$
This isomorphism also holds when $X$ is a stack and $X^\an$ is regarded as an
orbifold, \cite{oda,zoonekynd}. It follows that there are natural isomorphisms
$$
\pi_1^\prol(X\otimes \Kbar,x) \cong \pi_1^\prol(X^\an,x^\an) 
\text{ and }
\pi_1^\un(X\otimes \Kbar,x) \cong \pi_1^\un(X^\an,x^\an)_{/\Ql}. 
$$

When $G$ is a topological group and $V$ is a continuous $G$-module, 
$H^\dot(G,V)$ will denote the continuous cohomology (in the sense of Tate
\cite{tate}) of $G$ with coefficients in $V$.

\section{Monodromy Representations on Fundamental Groups}
\label{sec:arith}

\subsection{Monodromy representations}

Let $C \to B$ be a family of proper smooth pointed curves of genus $g$
defined over $K$. In our terminology, this means that $B$ is a stack, $C \to
B$ is proper and smooth and endowed with a section $\xi:B \to C$, where each
geometric fiber is a smooth, proper curve of genus $g$.

Assume that $B$ is geometrically connected. Let $x$ be a geometric point of
$B$, and $C_x$ the fiber over $x$. Suppose that $g \geq 1$. Then we have an
exact sequence of profinite groups
\begin{equation}
\label{eq:SGA1}
1 \to \pi_1^\alg(C_x,\xtilde) \to \pi_1^\alg(C,\xtilde) \to \pi_1^\alg(B,x)
\to 1,
\end{equation}
where $\xtilde = \xi(x)$.

The {\em monodromy representation}
\begin{equation}
\label{eq:monod}
\rho_{C,x}^\alg : \pi_1^\alg(B,x) \to \Aut\pi_1^\alg(C_x,\xtilde)
\end{equation}
on the fundamental group associated to the data $C\to B$, $\xi$ and $x$
is defined to be the composite
$$
\begin{CD}
\pi_1^\alg(B,x) @>{\xi_\ast}>> \pi_1^\alg(C,\xtilde) 
@>{\mathrm{conjugation}}>> \Aut \pi_1^\alg(C_x,\xtilde).
\end{CD}
$$

By functoriality, this homomorphism also induces continuous actions
\begin{equation}
\label{eq:monod:*}
\rho_{C,x}^\prol : \pi_1^\alg(B,x) \to \Aut \pi_1^\prol(C_x,\xtilde),\quad
\rho_{C,x}^\un : \pi_1^\alg(B,x) \to \Aut \pi_1^\un(C_x,\xtilde)
\end{equation}
of $\pi_1^\alg(B,x)$ on the pro-$\l$ and continuous unipotent
completions, respectively, of $\pi_1^\alg(C_x,\xtilde)$. In particular, when
$B$ is $\Spec K$, $\pi_1^\alg(B,x)$ is $G_K$, the absolute Galois group of $K$,
and these representations become the Galois representations
$$
\rho_{C,x}^\ast : G_K \to \Aut \pi_1^\ast(C_x,\xtilde)
$$
where $\ast \in \{\alg,\prol,\un\}$.

\subsection{The universal monodromy representation}
\label{sec:univ-monod}

Suppose that $g \geq 1$. Denote the moduli stack of pointed smooth projective
curves of genus $g>0$ over $K$ by $\M_g^1$. Denote the universal family over it
by $\cC_g^1 \to \M_g^1$, and its tautological section by $\xihat:\M_g^1 \to
\cC_g^1$. These are all defined over $\Spec K$.

Each pointed curve $(C,x)$ defined over an algebraically closed extension of
$K$ determines a geometric point $[C,x]$ of $\M_g^1$; the fiber of the pointed
universal curve over this point may be identified with $(C,x)$. Associated to
this geometric point is the monodromy representation
\begin{equation}
\label{eq:univ_monod}
\rhohat_{C,x}^\alg : \pi_1^\alg(\M_g^1,[C,x]) \to \Aut \pi_1^\alg(C,x).
\end{equation}
This is called the {\it universal monodromy representation} for the pointed
curve $(C,x)$ for reasons that we now recall.

For each proper family $\xymatrix@1{C \ar[r] & B\ar@/_/[l]|-{\xi}}$
of smooth pointed curves of genus $g$ defined over $K$, there is a unique
morphism $[C,\xi]: B \to \M_g^1$, called the {\em classifying morphism} of
the family, such that the pointed family
$$
\xymatrix{
C \ar[r] & B\ar@/_/[l]|-{\xi}
}
\quad\text{is the pullback of}\quad
\xymatrix{
\cC_g^1 \ar[r] & \M_g^1 \ar@/_/[l]|{\xihat}
}
$$
along $[C,\xi]$. It takes the geometric point $x$ of $B$ to $[C_x,\xtilde]$.
By universality, the fiber of $\cC_g^1 \to \M_g^1$ over $[C,\xi](x)$ can
be identified naturally with $(C_x,\xtilde)$. Thus, the diagram
$$
\xymatrix{
\pi_1^\alg(B,x)\ar[r]^(0.4){\rho_{C,x}^\alg}\ar[d]_{[C,\xi]_\ast} &
\Aut \pi_1^\alg(C_x,\xtilde)\cr
\pi_1^\alg(\M_g^1,[C_x,\xtilde]) \ar[ur]_{\rhohat_{C_x,\xtilde}^\alg}
}
$$
commutes, so that the monodromy representation of $C \to B$ associated to
the geometric point $x$ of $B$ factors through the universal monodromy
representation (\ref{eq:univ_monod}) associated to $(C_x,\xtilde)$.

By functoriality, the universal monodromy representation (\ref{eq:univ_monod})
also has pro-$\l$ and prounipotent incarnations:
\begin{equation}
\label{univ_monod:l}
\rhohat_{C,x}^\prol : \pi_1^\alg(\M_g^1,[C,x]) \to \Aut \pi_1^\prol(C,x)
\end{equation}
and
\begin{equation}
\label{univ_monod:un}
\rhohat_{C,x}^\un : \pi_1^\alg(\M_g^1,[C,x]) \to \Aut \pi_1^\un(C,x)
\end{equation}
through which the representations (\ref{eq:monod:*}) factor.

\subsection{Geometric monodromy}
Suppose that $(C,x)$ is a pointed curve defined over an algebraically closed
extension field of $K$. This determines a geometric point $[C,x]$ of $\M_g^1$.
It is standard that there is an exact sequence
\begin{equation}
\label{eq:ext-MgK}
1 \to \pi_1^\alg(\M_g^1\otimes \Kbar,[C,x]) \to \pi_1^\alg(\M_{g}^1,[C,x])
\to G_K \to 1. 
\end{equation}
The restriction
$$
\theta_{C,x}^\alg :
\pi_1^\alg(\M_g^1\otimes \Kbar,[C,x]) \to \Aut \pi_1^\alg(C,x)
$$
of the universal monodromy (\ref{eq:univ_monod}) to the geometric fundamental
group is called the {\it geometric monodromy} of the universal curve.

By functoriality, $\theta_{C,x}^\alg$ induces representations
$$
\theta_{C,x}^\prol :
\pi_1^\alg(\M_g^1\otimes \Kbar,[C,x]) \to \Aut \pi_1^\prol(C,x)
$$
and
$$
\theta_{C,x}^\un :
\pi_1^\alg(\M_g^1\otimes \Kbar,[C,x]) \to \Aut \pi_1^\un(C,x).
$$

Now suppose that $K$ is a subfield of $\C$ and that $(C,\xi)$ is a pointed
curve defined over $K$. Let $\Kbar$ be the algebraic closure of $K$ in $\C$. We
regard $\M_g^{1\,\an}$ as an orbifold. There is a natural isomorphism
$$
\G_\Canx \cong \pi_1(\M_g^{1\,\an},[\Canx]),
$$
which is unique up to conjugation by an element of $\Aut(\Canx)$. Consequently,
there is a natural isomorphism
$$
\Ghat_\Canx \cong
\pi_1^\alg(\M_g^1\otimes \Kbar,[\Cbar,\overline{\xi}]),
$$
with the profinite completion of $\G_\Canx$, where $\Cbar = C\otimes\Kbar$
and $\overline{\xi} = \xi\otimes \Kbar$. It is unique up to conjugation by an
element of $\Aut(\Canx)$.

By a result of Oda \cite{oda}, with respect to these identifications, the
geometric monodromy $\theta_{\Cbar,\overline{\xi}}$ is the completion of the
tautological action
$$
\G_{\Canx} \to \Aut \pi_1(\Canx).
$$

\begin{remark}
Suppose that $(\Cx)$ is a pointed curve over $K$. Let $x$ be a geometric point
of $\Spec K$. A fundamental problem (cf.\ \cite{MAPGAL}) is to understand the
image of the representation $\rho_{C,x}^\ast : G_K \to \Aut
\pi_1^\ast(C_x,\xtilde)$, where $\ast \in \{\alg,\prol,\un\}$. Since the
diagram
$$
\xymatrix{
G_K \ar[r]^(0.3){[\Cx]_\ast}\ar[dr]_{\rho_{C,x}^\ast} &
\pi_1^\alg(\M_g^1,[C_x,\xtilde])\ar[d]^{\rhohat_{C_x,\xtilde}^\ast}\cr
& \Aut \pi_1^\ast(C_x,\xtilde).
}
$$
commutes, we have the inclusions
$$
\im \rho_{C,x}^\ast \hookrightarrow \im \rhohat_{C_x,\xtilde}^\ast
\hookrightarrow \Aut \pi_1^\ast(C_x,\xtilde).
$$
One can ask if the left hand inclusion can be an isomorphism. In \cite{MAPGAL},
the second author and Tamagawa show that when $C_x$ is an {\it affine} curve of
$g\neq 1$ and $*=\alg$, this inclusion is never an equality. More precisely,
they prove that when $K$ is a number field, no nontrivial element of the
geometric part $\theta_{C_x,\xtilde}^\alg$ lies in the image of
$\rho_{C,x}^\alg$. On the other hand, when $\ast = \prol$, equality holds for
the general curve (in the sense of Hilbert's Irreducibility Theorem).
Theorem~\ref{th:main-aut} treats the case $*=\un$, where image is replaced by
the Zariski closure of the image.
\end{remark}

\section{Cycle Classes}
\label{sec:univ_classes}

In this section, we describe the class of the cycle class in the relative
jacobian associated to a pointed family of curves $C \to B$.

\subsection{The $\l$-adic Abel-Jacobi map}
\label{sec:Abel-Jacobi}

Let $B$ be a geometrically connected stack over $K$, and let $f: X \to B$ be a
proper smooth morphism locally of finite type. Let $\alpha$ be a relative
algebraic cycle in $X/B$ of codimension $r$. That is, $\alpha$ is a formal
integral linear combination of codimension-$r$ closed subvarieties of $X$
equidimensional over $B$. Denote the group of such relative cycles by
$\Cyc^r(X/B)$.

One has the cycle class map
$$
cl : \Cyc^r(X/B) \to \Het^0(B,R^{2r}f_*\,\Z_l(r)).
$$
When $K$ is a subfield of $\C$, we also have its analogue
$$
cl^\an : \Cyc^r(X^\an/B^\an) \to H^0(B^\an,R^{2r}f^\an_*\,\Z)
$$
for the associated analytic family $f^\an : X^\an \to B^\an$ by
$\Cyc^r(X^\an/B^\an)$.

The kernels of these cycle maps are, by definition, the groups  by
$\Cyc^r_\hom(X/B)$ and $\Cyc^r_\hom(X^\an/B^\an)$, respectively, of
homologically trivial cycles. Elements of the first group may be detected via
the ($\l$-adic) Abel-Jacobi mapping
\begin{equation}
\label{eq:AJ}
A : \Cyc^r_\hom(X/B) \to \Het^1(B, R^{2r-1}f_*\,\Z_l(r)),
\end{equation}
and of the second by its topological analogue
$$
A^\an : \Cyc^r_\hom(X/B) \to H^1(B^\an, R^{2r-1}f^\an_*\,\Z).
$$
We now recall the definition of $A$.

Denote the support of an element $Z$ of $\Cyc^r_\hom(X/B)$ by $|Z|$, and its
inclusion into $X$ by $\iota: |Z| \hookrightarrow X$.  Define
$\cH^{2r}_{|Z|}(X,\Zl(r))$ to be the sheaf $(R^{2r}\iota^!)\Zl(r)$ over $B$, so
that  $\big({f|_{|Z|}}\big)_\ast\cH^{2r}_{|Z|}(X,\Zl(r))$ is the sheaf whose
stalk at $b\in B$ is $H^{2r}_{|Z|_b}(X_b, \Zl(r))$. Pulling back the Gysin
sequence
(cf.\
\cite{milne})
\begin{multline*}
0 \to R^{2r-1}f_*\,\Zl(r) \to R^{2r-1}{\big(f|_{X-|Z|}\big)}_*\,\Zl(r) \cr
\to \big({f|_{|Z|}}\big)_\ast\cH^{2r}_{|Z|}(X,\Zl(r))
\to R^{2r}f_*\,\Zl(r)
\end{multline*}
along the mapping $\Zl \to {f|_{|Z|}}_*\,\cH^{2r}_{|Z|}(X,\Zl(r))$ gives a
short exact sequence 
\begin{equation}
\label{eq:c-c-extension}
0 \to R^{2r-1}f_*\,\Zl(r) \to \cE \to \Zl \to 0.
\end{equation}
The class of this extension of smooth $\Zl$-sheaves over $B$ in
$$
\Ext^1(\Zl, R^{2r-1}f_*\,\Zl(r))\cong \Het^1(B, R^{2r-1}f_*\,\Zl(r))
$$
is $A(Z)$.

In the case where $B$ is the spectrum of the field $K$, the Abel-Jacobi map is
constant on rational equivalence classes of cycles, and therefore defines a
map
$$
A : CH^r_\hom(X/K) \to 
\Het^1(K, R^{2r-1}f_*\,\Zl(r))\cong H^1(G_K, H^{2r-1}(X,\Zl(r)))
$$
which coincides with the usual $\l$-adic Abel-Jacobi map (cf.\
\cite[Lem.~9.4]{jannsen}).

When $K \subseteq \C$, the topological Abel-Jacobi mapping is defined in an
analogous way. Standard comparison theorems imply that there is a commutative
diagram
$$
\begin{CD}
\Cyc^r(\Xbar/\Bbar) @>{cl}>> \Het^0(\Bbar,R^{2r}\bar{f}_*\,\Z_l(r)) \cr
@VVV @VV{\simeq}V \cr
\Cyc^r(X^\an/B^\an) @>{cl^\an}>> H^0(B^\an,R^{2r}f^\an_*\,\Z(r))\otimes\Zl
\end{CD}
$$
where $\bar{f} : \Xbar \to \Bbar$ is the pullback of $f : X \to B$
along $K \to \Kbar$. Once again, standard comparison theorems imply the
commutativity of the diagram
\begin{equation}
\label{compat}
\begin{CD}
\Cyc^r_\hom(\Xbar/\Bbar) @>{A}>> \Het^1(\Bbar,R^{2r-1}\bar{f}_*\,\Z_l(r)) \cr
@VVV @VV{\simeq}V \cr
\Cyc^r_\hom(X^\an/B^\an) @>{A^\an}>>
H^1(B^\an,R^{2r-1}f^\an_*\,\Z(r))\otimes\Zl.
\end{CD}
\end{equation}

\subsection{Setting and notation}
\label{sec:setting}

Here we make precise the setting and notation we will use in much of the
rest of the paper.

Suppose that $C \to B$ is a proper, geometrically connected, smooth family of
curves of genus $g$ defined over $K$, where $g\ge 3$. Suppose that $\xi:B \to
C$ is a section. Denote the relative jacobian $\Jac_{C/B}$ by $f: X \to B$.
Suppose that $x$ is a geometric point of $B$.

For all such families, we shall denote the $\l$-adic local system
$R^{2g-3}f_\ast \Zl(g-1)$ over $B$ by $\L_\Zl$. Similarly, denote the $\l$-adic
local system $R^{2g-1}f_\ast \Zl(g)$ over $B$ by $\H_\Zl$. Note that $\H_\Zl$
is the local system corresponding to the action of $\pi_1^\alg(B,x)$ on the
module $H_\Zl := \Het^1(C_x,\Zl(1))$ and that $\L_\Zl$ is the local system
corresponding to the action of $\pi_1^\alg(B,x)$ on the module $L_\Zl :=
\big[\Lambda^3 H_\Zl\big](-1)$.

If $B$ is geometrically connected and $x$ is a geometric point of $B$, then
there are natural isomorphisms (cf.\ \cite{milne})
$$
\Het^1(B, \L_\Zl) \cong H^1(\pi_1^\alg(B,x), L_\Zl) \text{ and }
\Het^1(B, \H_\Zl) \cong H^1(\pi_1^\alg(B,x), H_\Zl).
$$

When $K$ is a subfield of $\C$, we can associate to the corresponding analytic
families $C^\an \to B^\an$ and $f^\an : X^\an \to B^\an$ the local systems
$$
R^{2g-3}f^\an_\ast \Z(g-1) \text{ and } R^{2g-1}f^\an_\ast \Z(g),
$$
which we shall denote by $\L_\Z$ and $\H_\Z$, respectively.

\subsection{The class $\mu(\Cx)$}
\label{sec:C-C-}

As in Section~\ref{sec:introduction}, we have two closed embeddings of $C$ in
its Jacobian $X$; namely,  $y \mapsto [y]-[\xi]$ and $y \mapsto [\xi]-[y]$.
These define the two relative algebraic cycles
$$
C_\xi,\, C_\xi^- \in \Cyc^{g-1}(X/B),
$$
respectively. Denote their difference $C_\xi - C_\xi^-$ by $\alpha(C,\xi)$.
Since the endomorphism $-\id$ of $X$ multiplies $\alpha(C,\xi)$ by $-1$ and
acts trivially on $R^2f_*\Zl$, which is torsion free, it follows that
$\alpha(C,\xi)$ is homologically trivial. It therefore determines a class
$$
A(\alpha(C,\xi)) \in \Het^1(B, \L_\Zl)
$$
which we shall denote by $\mu(C,\xi)$.

We will consider $\mu(C,\xi)$ as an element of $H^1(\pi_1^\alg(B,x), L_\Zl)$.
This generalizes the definition given in \S\ref{sec:introduction}.

This class is functorial with respect to morphisms of pointed families. In
particular, we have the class
$$
\mu(\cC_g^1,\xihat) \in \Het^1(\M_g^1,\L_\Zl)
$$
of the universal curve, which is universal.

\begin{proposition}
\label{prop:C-specialization}
If $\xymatrix@1{C \ar[r] & B\ar@/_/[l]|-{\xi}}$ is a pointed family over $K$
and $x$ is a geometric point of $B$, then
$$
\mu(C,\xi) \in H^1(\pi_1^\alg(B,x),L_\Zl)
$$
is the pullback of
$$
\mu(\cC_g^1,\xihat) \in H^1(\pi_1^\alg(\M_g^1,[C_x,\xtilde]),L_\Zl)
$$
under the group homomorphism induced by the classifying morphism
$$
[C,\xi] : (B,x) \to (\M_g^1,[C_x,\xtilde]). \qed
$$
\end{proposition}

When $K$ is a subfield of $\C$, we have the relative cycle $C^\an_\xian -
C^{\an-}_\xian$ in $f^\an: X^\an \to B^\an$. This determines a class
$$
\mu^\an(C,\xi) \in H^1(B^\an,\L_\Z).
$$

Applying this construction to the universal pointed curve $\cC_g^{\an\,1} \to
\M_g^{\an\,1}$, we obtain the universal such class
$$
\mu^\an(\cC_g^1,\xihat) \in H^1(\M_g^{\an\, 1},\L_\Z).
$$

Since the diagram (\ref{compat}) commutes and since this class is respected by
base change, we have:

\begin{proposition}
\label{prop:C-comparison}
The class $\mu^\an(C,\xi)$ corresponds to $\mu(C,\xi)$ under the comparison
isomorphism
$$
\Het^1(B\otimes \Kbar,\L_\Zl) \cong \Het^1(B\otimes \C,\L_\Zl)
\cong H^1(B^\an,\L_\Z)\otimes \Zl.
$$
Moreover, the class $\mu^\an(C,\xi)$ is the pullback of
$\mu^\an(\cC_g^{\an\,1},\xihat)$ along the classifying mapping $[C,\xi] :
B^\an\to \M_g^{\an\,1}$. \qed
\end{proposition}

\subsection{The class $\kappa(\Cx)$}

Denote the canonical divisor class of $C$ over $B$ by $\K_C$. The cycle
$(2g-2)\xi - \K_C$ is homologous to zero on each geometric fiber of $C \to B$.
It therefore defines a class
$$
\kappa(\Cx) \in \Het^1(B,\H_\Zl).
$$

Applying this construction to the universal pointed curve $\cC_g^1 \to \M_g^1$,
we obtain the universal such class
$$
\kappa(\cC_g^1,\xihat) \in \Het^1(\M_g^1,\H_\Zl).
$$
As in the case of $\mu(\Cx)$, this class pulls back to $\kappa(\Cx)$ along
the classifying morphism of the pointed family $C \to B$.

\begin{proposition}
\label{prop:K-specialization}
If $\xymatrix@1{C \ar[r] & B\ar@/_/[l]|-{\xi}}$ is a pointed family over $K$
and $x$ is a geometric point of $B$, then
$$
\k(C,\xi) \in H^1(\pi_1^\alg(B,x),H_\Zl)
$$
is the pullback of
$$
\k(\cC_g^1,\xihat) \in H^1(\pi_1^\alg(\M_g^1,[C_x,\xtilde]),H_\Zl)
$$
under the group homomorphism induced by the classifying morphism
$$
[C,\xi] : (B,x) \to (\M_g^1,[C_x,\xtilde]). \qed
$$
\end{proposition}

When $K$ is a subfield of $\C$, we have the relative cycle $(2g-2)\xian -
\K_C$ in $C \to B$. This determines a class
$$
\kappa^\an(C,\xi) \in H^1(B^\an,\H_\Z).
$$

Applying this construction to the universal pointed curve $\cC_g^{\an\,1} \to
\M_g^{\an\,1}$, we obtain the universal such class
$$
\kappa^\an(\cC_g^1,\xihat) \in H^1(\M_g^{\an\, 1},\H_\Z).
$$

The commutativity of the diagram (\ref{compat}) and naturality of the class
under base change implies:

\begin{proposition}
\label{prop:K-comparison}
The class $\kappa^\an(C,\xi)$ corresponds to
$\kappa(C,\xi)$ under the comparison isomorphism
$$
\Het^1(B\otimes \Kbar,\H_\Zl) \cong \Het^1(B\otimes \C,\H_\Zl)
\cong H^1(B^\an,\H_\Z)\otimes \Zl.
$$
Moreover, the class $\k^\an(C,\xi)$ is the pullback of
$\k^\an(\cC_g^{\an\,1},\xihat)$ along the classifying mapping $[C,\xi] :
B^\an\to \M_g^{\an\,1}$.
\end{proposition}

\subsection{The classes $\nu(C)$ and $\nutilde(\Cx)$}

Denote by $\nutilde(C,\xi)$ the image of $\mu(C,\xi)$ under the quotient
mapping
$$
\Het^1(B,\L_\Z) \to \Het^1(B,\L_\Z/\H_\Zl).
$$

\begin{proposition}
\label{prop:N-specialization}
If $\xymatrix@1{C \ar[r] & B\ar@/_/[l]|-{\xi}}$ is a pointed family over $K$
and $x$ is a geometric point of $B$, then
$$
\nutilde(C,\xi) \in H^1(\pi_1^\alg(B,x),L_\Zl/H_\Zl)
$$
is the pullback of
$$
\nutilde(\cC_g^1,\xihat) \in H^1(\pi_1^\alg(\M_g^1,[C_x,\xtilde]),L_\Zl/H_\Zl)
$$
under the group homomorphism induced by the classifying morphism
$$
[C,\xi] : (B,x) \to (\M_g^1,[C_x,\xtilde]). \qed
$$
\end{proposition}

When $K$ is a subfield of $\C$, define
$$
\nutilde^\an(C,\xi) \in H^1(B^\an,\L_\Zl/\H_\Zl),
$$
to be the image of $\mu^\an(C,\xi)$ under the natural mapping
$$
H^1(B^\an,\L_\Z) \to H^1(B^\an,\L_\Zl/\H_\Zl).
$$

Applying this construction to the universal pointed curve $\cC_g^{\an\,1} \to
\M_g^{\an\,1}$, we obtain the universal such class
$$
\nutilde^\an(\cC_g^1,\xihat) \in H^1(\M_g^{\an\, 1},\L_\Z/\H_\Z).
$$

Since the diagram (\ref{compat}) commutes and since the this class is respected
by base change, we have:

\begin{proposition}
\label{prop:N-comparison}
The class $\nutilde^\an(C,\xi)$ corresponds to $\nutilde(C,\xi)$ under the
comparison isomorphism
$$
\Het^1(B\otimes \Kbar,\L_\Zl/\H_\Zl) \cong \Het^1(B\otimes \C,\L_\Zl/\H_\Zl)
\cong H^1(B^\an,\L_\Z/\H_\Z)\otimes \Zl.
$$
Moreover, the class $\nutilde^\an(C,\xi)$ is the pullback of
$\nutilde^\an(\cC_g^{\an\,1},\xihat)$ along the classifying mapping $[C,\xi] :
B^\an\to \M_g^{\an\,1}$. \qed
\end{proposition}

Although we have used the rational point $\xi$ to define the class
$\nutilde(C,\xi)$, this class depends only on $C$ and not on the rational point
$\xi$ (or even the existence of the rational point). The following result is
proved in Section~\ref{sec:nu}.

\begin{proposition}
\label{prop:equal}
The class $\nutilde(C,\xi)$ depends only on $C\to B$ and is independent of the
choice or existence of the $B$-rational point $\xi$.
\end{proposition}

For this reason, we shall denote $\nutilde(C,\xi)$ by $\nu(C)$. Note, however,
that our proof of Theorem~\ref{th:main-aut} does not depend on this equality,
and that Section~\ref{sec:nu}, where this equality is proved, can be read
at this point of the paper.

\section{The Universal Johnson Class}
\label{sec:univ_johnson}

In this section we construct a universal cohomology class
$$
\mhat^\prol \in H^1(\Aut \pi_1^\prol(S,x),L_\Zl)
$$
where $S$ is a compact Riemann surface of genus $g \ge 3$ and $L_\Zl$ is the
third exterior power of $H_1(S,\Zl)$. This construction is a cohomological
interpretation of an $\l$-adic version of the Johnson homomorphism
\cite{johnson:homom}, which is, in turn, an outgrowth of a construction of
Magnus.

The class $\mhat^\prol$ is a characteristic class for groups $G$ that act
continuously on $\pi_1^\prol(S,x)$, such as Galois groups. Indeed, the class
$\mhat^\prol$ can be pulled back along the homomorphism
$$
G \to \Aut \pi_1^\prol(S,x)
$$
that defines the action to give a class in $H^1(G,L_\Zl)$.

\subsection{The Magnus homomorphism}

Suppose that $\pi$ is a topological group. The examples we have in mind are
where $\pi$ is either a discrete group (viewed as a topological group with the
discrete topology) or a profinite group. Denote its group of continuous
automorphisms by $\Aut \pi$. We shall denote the lower central series
filtration of $\pi$ by
$$
\pi = L^1 \pi \supseteq L^2 \pi \supseteq L^3 \pi \supseteq \cdots
$$
where $L^{k+1}\pi$ is the closure of $[L,L^k]$.

Set $H_1(\pi) = \pi/L^2\pi$. Denote the kernel and image of the natural
homomorphism
$$
\Aut \pi \to \Aut H_1(\pi)
$$
by $T$ and $R$, respectively, where $H_1(\pi)$ denotes the Hausdorff
abelianization of $\pi$, which is the quotient of $\pi$ by the closure of its
commutator subgroup. Since the filtration $L^\dot$ of $\pi$ is by
characteristic subgroups, there is a natural homomorphism
$$
\etilde : T \to \Hom(H_1(\pi),\Gr^2_L \pi)
$$
into the group of continuous homomorphisms $H_1(\pi) \to \Gr^2_L \pi$. It is
defined by taking the element $\phi$ of $T$ to the homomorphism
$$
u \mapsto \phi(\utilde)u^{-1} \mod L^3\pi
$$
where $\utilde \in \pi$ is any lift of $u \in H_1(\pi)$. Since $\phi$ is
continuous and acts trivially on $H_1(\pi)$, the homomorphism $\etilde(u)$ is
easily seen to be well defined and continuous. It is also easy to see that
$\etilde$ is itself a continuous homomorphism.

Each graded quotient of the lower central series of $\pi$ is naturally a
continuous $R$-module. It is an exercise to show that the homomorphism
$$
\e_1 : H_1(T) \to \Hom(H_1(\pi),\Gr^2_L \pi)
$$
induced by $\etilde$ is continuous and $R$-linear. We shall call it the {\it
Magnus homomorphism} as Magnus \cite{magnus} was the first to study the kernel
of the homomorphism from the Automorphism group of a free group to the
automorphisms of its abelianization. A similar construction was used by
Andreadakis \cite{andreadakis} when studying the automorphism group of a free.
We shall view it as an element
$$
\e_1 \in \Hom_R(H_1(T),\Hom(H_1(\pi),\Gr^2_L \pi)),
$$
which we shall call the {\it Magnus element}.

It is well known that $\Gr_L^\dot \pi$ is a graded Lie algebra (cf.\
\cite{serre:LALG}).  In our case, the bracket
$$
[\blank,\blank] : \Gr_L^n \pi \otimes \Gr_L^m \pi \to \Gr_L^{n+m}\pi
$$
is easily seen to be $R$-equivariant.

The definition of $\e_1$ can be extended to define continuous $R$-equivariant
mappings
$$
\e_n : H_1(T) \to \Hom(\Gr_L^n \pi,\Gr_L^{n+1}\pi)
$$
in such a way that, for each $\phi \in H_1(T)$, $(\e_n(\phi))_n \in \End
\Gr_L^\dot \pi$ is a derivation in the following sense.

\begin{proposition}
\label{derivation}
If $u \in \Gr_L^i \pi$, $v \in \Gr_L^j \pi$ and $i+j = n$, then for each
$\phi \in H_1(T)$
$$
\e_n(\phi)([u,v]) = [\e_i(\phi)(u),v] + [u,\e_j(\phi)(v)]. \qed
$$
\end{proposition}

\subsection{The Johnson homomorphism}
\label{sec:johnson}

Now suppose that $\pi$ is either a hyperbolic surface group (that is, the
fundamental group of a compact orientable surface of genus $g \ge 2$) with the
discrete topology, or the pro-$\l$ completion of such a surface group with the
pro-$\l$ topology. Let $A$ denote $\Z$ in the first case, and $\Zl$ in the
second. Then $H_A := H_1(\pi)$ is a free $A$-module of rank $2g$, which is
endowed with a unimodular skew symmetric bilinear form $\w$ --- the
intersection form. In this case, $R = \GSp(H_A)$, the group of symplectic
similitudes of $H_A$. It is an extension
$$
\begin{CD}
1 @>>> \Sp(H_A) @>>> \GSp(H_A) @>{\chi}>> A^\times @>>> 1
\end{CD}
$$
of the units of $A$ by the symplectic group $\Sp(H_A)$. The character $\chi$
is the 1-dimensional representation of $\GSp$ given by its action on the
symplectic form $\w$.

By a theorem of Labute \cite{labute}, the Lie algebra $\Gr_L^\dot \pi$ is
torsion free as an $A$-module, and the natural surjection
$$
\Lib(H_A) \to \Gr_L^\dot \pi
$$
induces a Lie algebra isomorphism
$$
\Lib(H_A)/(q) \cong \Gr_L^\dot \pi
$$
induced by the isomorphism $H_A \to \Gr_L^1 \pi$, where $\Lib(V)$ denotes the
free Lie algebra over $A$ generated by the $A$-module $V$, and $q$ is the
symplectic element of $\Lambda^2 H_A$. In particular, the bracket $H_A\otimes
H_A \to \Gr_L^2 \pi$ is surjective and induces an $\GSp(H_A)$-equivariant
isomorphism
$$
\Gr_L^2 \pi \cong \big(\Lambda^2 H_1(\pi)\big)/q.
$$

Denote the field of fractions of $A$ by $F$. That is $F = \Q$ when  $A = \Z$
and $F = \Ql$ when $A = \Zl$. Then $H_F := H_A\otimes F$ is the fundamental
representation of $\GSp(H_F)$. Since the bracket
$$
H_F \otimes \big[\big(\Gr_L^n \pi\big) \otimes F\big]
\to \big(\Gr_L^{n+1} \pi\big) \otimes F
$$
is surjective and $\GSp(H_F)$-equivariant for all $n\ge 1$, it follows by
induction on $n$ that each
$$
\big(\Gr_L^n \pi\big) \otimes F
$$
is also a $\GSp(H_F)$-module.

Henceforth we shall regard $\GSp(H)$ as an algebraic group over $\Z$. Its group
of $R$-rational points being $\GSp(H\otimes_\Z R)$. The action of this group on
the symplectic form defines a homomorphism $\GSp(H) \to \Gm$. Denote the
pullback of the standard representation of $\Gm$ by $\Z(1)$. Define $\Z(n)$ to
be its $n$th tensor power when $n > 0$, to be trivial when $n=0$,  and to be
the dual of $\Z(-n)$ when $n < 0$. The tensor product of a representation $V$
of $\GSp(H)$ with $\Z(n)$ will be denoted by $V(n)$. (Later, this will
correspond to Tate twist.) Since the pairing $H_A \otimes H_A \to A(1)$ is
$\GSp$-equivariant, so is its adjoint, the isomorphism
$$
H_A(-1) \cong \Hom(H_A,A).
$$

In this case, the Magnus element induces a homomorphism
\begin{equation}
\label{eq:magnus}
\tautilde : H_1(T) \to H_A(-1) \otimes \big[\big(\Lambda^2 H_A\big)/q\big].
\end{equation}

Johnson made the further observation that when $\pi$ is a surface group, the
Magnus homomorphism has restricted image. First note that there is a natural
$\GSp(H_A)$-equivariant embedding
$$
f : \big[\Lambda^3 H_A\big](-1) \to
H_A(-1) \otimes \big[\big(\Lambda^2 H_A\big)/q\big]
$$
defined by
$$
a \wedge b \wedge c \mapsto
a\otimes (b\wedge c) + b\otimes (c\wedge a) + c\otimes (a\wedge b).
$$

\begin{lemma}[Johnson]
\label{lem:factor}
The image of $\tautilde$ lies in the image of $f$.
\end{lemma}

Because of this, the Magnus homomorphism $\tautilde$ induces an
$\GSp(H_A)$-equivariant homomorphism
$$
\tau : H_1(T) \to \Lambda^3 H_A(-1).
$$
This is the celebrated {\it Johnson homomorphism}.

Since Johnson's Lemma is novel in the $\l$-adic case, where $T$ is much larger
than the Torelli group, and because we are applying it in such situations, we
will give a proof. First we need:

\begin{lemma}
\label{lem:jacobi-low}
The sequence of $\GSp(H_A)$-modules
$$
\begin{CD}
0 @>>>  \Lambda^3 H_A @>f(1)>> H_A \otimes \Gr_L^2 \pi
@>{\mathrm{bracket}}>> \Gr_L^3 \pi @>>> 0
\end{CD}
$$
is exact.
\end{lemma}

\begin{proof}
First, after tensoring the sequence with $F$, it is exact. This is well known
(cf.\ \cite[8.3]{hain:torelli}) and easily proved using the representation
theory of $\Sp(H_F)$. The right hand term is irreducible, and the left hand
term is the sum of two irreducible components when $g \ge 3$ and is isomorphic
to $H_F$ when $g=2$.

The result now follows over $A$ as, by Labute, each $\Gr_L^n \pi$ is torsion
free, since the right hand mapping is surjective by the definition of
the filtration $L^\dot$, and since the kernel of the bracket divided by
the image of $f$ is torsion free, which is easily checked by constructing
an $A$-module splitting, for example.
\end{proof}

\begin{proof}[Proof of Lemma~\ref{lem:factor}]
We use the fact that $H_1(T)$ acts on $\Gr_L^\dot \pi$ as derivations.
By Proposition~\ref{derivation}, each $\phi \in T$ induces a map of short
exact sequences
$$
\xymatrix{
0 \ar[r] & A\cdot q \ar[r] \ar@{.>}[d] & \Lambda^2 H_A
	\ar[r]^{\mathrm{bracket}} \ar[d]^{\sigma(\phi)} & \Gr_L^2 \pi \ar[r]
	\ar[d]^{\e_2(\phi)} & 0 \cr
0 \ar[r] & \Lambda^3 H_A \ar[r]^(.4){f(1)} & H_A \otimes \Gr_L^2 \pi
	\ar[r]^(.6){\mathrm{bracket}} & \Gr_L^3 \pi \ar[r] & 0 \cr
}
$$
where $\sigma(\phi)(a\wedge b) = a \otimes \e_1(\phi)(b) - b \otimes
\e_1(\phi)(a)$. Since the function $\e_1(\phi) : H_A \to \Gr^2_L \pi$
corresponds, by duality, to
$$
\tautilde(\phi)=
\sum_{j=1}^g \big(a_j\otimes \e_1(\phi)(b_j) - b_j\otimes \e_1(\phi)(a_j)\big)
\in H_A\otimes \Gr^2_L \pi,
$$
we see that $\tautilde(\phi)= \sigma(\phi)(q)$. The definitions then imply that
\begin{multline*}
\mathrm{bracket} \circ \sigma(\phi)(q) = \mathrm{bracket} \circ
\sigma(\phi)\bigg(\sum_{j=1}^g a_j \wedge b_j\bigg) \cr
= \sum_{j=1}^g \big([\e_1(\phi)(a_j),b_j] - [a_j, \e_1(\phi)(b_j)]\big)
= \e_2(\phi)\bigg(\sum_{j=1}^g [a_j,b_j]\bigg) = 0.
\end{multline*}
\end{proof}

\subsection{A characteristic class for groups acting on surface groups}
\label{sec:char_class}

We retain the notation of the preceding paragraphs. In this paragraph, we
explain how the Johnson homomorphism, or more accurately, twice the Johnson
homomorphism, can be described as a cohomology class in
$$
H^1(\Aut\pi, L_A),
$$
where $L_A = (\Lambda^3 H_A)(-1)$. This can be viewed as a universal
characteristic class of the group action.

One would like to say that the extension
$$
1 \to T \to \Aut \pi \to \GSp(H_A) \to 1
$$
gives rise to a spectral sequence
$$
H^s(\GSp(H_A),H^t(T,L_A)) \implies H^{s+t}(\Aut \pi, L_A).
$$
This is true in the discrete case, but may not be in the profinite case (cf.\
\cite{jannsen:cont}). Nonetheless, one can easily establish, by a cocycle
computation, the exactness of the sequence
\begin{multline}
\label{eq:first-three}
\begin{CD}
0 @>>> H^1(\GSp(H_A),L_A) @>>> H^1(\Aut \pi,L_A)
\end{CD}
\cr
\begin{CD}
@>>> H^0(\GSp(H_A),H^1(T,L_A)) @>{d_2}>> H^2(\GSp(H_A),L_A)
\end{CD}
\end{multline}
which one would have should the spectral sequence exist.

Note that the $\GSp(H_A)$-modules $L_A$, $H_A$ and $L_A/H_A$ all have the
property that $-I \in \GSp(H_A)$ acts as $-\id$ on them.

\begin{lemma}
\label{lem:gen-vanish}
Suppose that $R_A$ is either $\GSp(H_A)$ or $\Sp(H_A)$ and that $g\ge 1$. If
$V$ is an $R_A$-module on which $-I\in R_A$ acts as $-\id$, then the groups
$H^k(R_A,V)$ have exponent $2$ (that is, they are annihilated by $2$). In
particular, these groups vanish if $A = \Zl$ when $\l \neq 2$. If $V$ is
torsion free, then there is a natural isomorphism
$$
H^1(R_A,V) \cong H^0(R_A,V/2).
$$
In particular, $H^1(R_A,L_A)$ vanishes for all $g\ge 1$.
\end{lemma}

\begin{proof}
The element $-I$ of the center of $R_A$ acts as $-\id$ on $H^k(R_A,L_A)$. But,
by ``center kills,'' elements of the center of $R_A$ act trivially on all
cohomology groups of $R_A$. It follows that $H^\dot(R_A,V)$ is annihilated by
$2$. Since $2$ is a unit of $\Zl$ when $\l \neq 2$, this implies the vanishing
of $H^\dot(R_A,L_A)$ in that case. It also implies the vanishing of
$H^0(R_A,V)$ when $V$ is torsion free.

If $V$ is torsion free, the Bockstein sequence associated to
$$
0 \to V \stackrel{2}{\to} V \to V/2 \to 0
$$
yields a natural isomorphism $\beta : H^0(R_A,V/2) \to H^1(R_A,V)$.

When $g=1$, $L_A = 0$, so that $H^1(R_A,L_A)$ vanishes. When $g=2$, cupping
with the symplectic form gives a natural isomorphism $H_A \cong L_A$. One can
easily verify that $H_A/2$ has no invariants. Indeed, fix a curve $C$ of genus
2. It is hyperelliptic. There are $2^4 - 1 = 15$ non-zero points of order 2 in
the Jacobian of $C$. These are the 15 differences of  2 distinct Weierstrass
points. Since none are fixed by the mapping class group, there are no
invariants. The case $g\ge 3$ is treated in \cite[\S6]{hain-reed}. There it is
proved that, when $A=\Z$, $H^0(R_A,L_A/2)$ vanishes. Since $L_A$ is torsion
free and $H^0(R_A,L_A/2)$ vanishes for all $g\ge 1$, $H^1(R_A,L_A)$ also
vanishes.
\end{proof}

The following generalizes some results of Morita \cite{morita} from mapping
class groups to the larger group $\Aut \pi$.

\begin{proposition} 
\label{prop:univ_class}
There is a unique class $\mhat \in H^1(\Aut\pi, L_A)$ whose image in
$H^0(\GSp(H_A),\Hom(H_1(T),L_A))$ is $2\tau$.
\end{proposition}

\begin{proof}
This follows from the facts that $H^1(\GSp(H_A),L_A)$ vanishes
(Prop.~\ref{prop:vanish}) and that $H^2(\GSp(H_A),L_A)$ has exponent $2$
(Lem.~\ref{lem:gen-vanish}) using the exact sequence (\ref{eq:first-three}).
\end{proof}

\begin{remark}
The image of $\tau$ under $d_2 : H^1(\Aut \pi,L_A) \to H^2(\GSp(H_A),A)$ is
non-zero when $2$ is not a unit in $A$ (i.e., when $A=\Z_2$ and $\Z$). To prove
this, it suffices to consider the case where $\pi$ is discrete and $A=\Z$ as
there is a homomorphism $\Aut \pi \to \Aut \pi^{(2)}$. In
\cite[p.~220]{morita}, Morita proves the equivalent statement that $\tau$ is
not in the image of
$$
H^1(\Aut \pi, L_\Z) \to \Hom(H_1(T), L_\Z)
$$
but that $2\tau$ is.
\end{remark}

We shall adopt the following notation for the universal class in the $\l$-adic
case:
\begin{equation}
\label{eq:def-of-ml}
\mhat^\prol \in H^1(\Aut \pi^\prol, L_\Zl).
\end{equation}

\begin{definition}
Suppose that $(C,\xi)$ is a pointed curve over $K$ of genus $g\ge 2$, where
$K$ is a subfield of $\C$. Define the class
$$
m(C,x) \in H^1(G_K,L_\Zl)
$$
to be the pullback of the universal class $\mhat^\prol$ along the homomorphism
$$
\rho_{C,x}^\prol : G_K \to \Aut \pi_1^\prol(\Canx).
$$
\end{definition}

This is the invariant of the Galois action on the fundamental group that
appears in the statement of the $\l$-adic Harris-Pulte Theorem in the
Introduction.

\subsection{Further torsion computations}
\label{sec:torsion-comps}
In this section Lemma~\ref{lem:gen-vanish} is extended to obtain more torsion
vanishing examples that are needed in subsequent sections.

Wedging with the symplectic element $q$ of $H_A$ defines a
$\GSp(H_A)$-invariant imbedding
$$
i : H_A \hookrightarrow L_A.
$$
There is also a $\GSp(H_A)$-invariant projection $c: L_A \to H_A$ defined by
$$
c : x \wedge y \wedge z \mapsto \w(x,y)z + \w(y,z)x + \w(z,x)y
$$
It is easy to verify that $c\circ i = (g-1)\id$ so that, when $g-1$ is a
unit in $A$, we have a $\GSp(H_A)$-module decomposition
\begin{equation}
\label{eq:decomp}
(c,p) : L_A \stackrel{\simeq}{\longrightarrow} H_A \oplus L_A/H_A,
\end{equation}
where $p$ denotes the projection $L_A \to L_A/H_A$.

When $A$ is a field of characteristic zero, $H_A$ and $L_A/H_A$ are irreducible
representations of the algebraic group $\Sp(H_A)$.

\begin{proposition}
\label{prop:vanish}
If $R_A$ is $\GSp(H_A)$ or $\Sp(H_A)$, then the groups $H^1(R_A,H_A)$,
$H^1(R_A,L_A)$ and $H^1(R_A,L_A/H_A)$ all vanish for all $g\ge 1$.
\end{proposition}

\begin{proof}
In view of Lemma~\ref{lem:gen-vanish}, we need only show that $H^0(R_A,V/2)$
vanishes when $V=H_A$ and $L_A/H_A$. It is easy to check (and well known) that
$H^0(R_A,H_A/2)$ vanishes for all $g\ge 1$. Indeed, when $g=1$ and $A = \Z$ or
$\Z_2$, then the image of $R_A$ in $\Aut (H_A/2)$ is $SL_2(\F_2)$ and $H_A/2$
is the two dimensional vector space over $\F_2$. The result follows in this
case. When $g > 1$, one can restrict to the subgroup $R'$ of $R_A$ that is a
product of $g$ copies of $SL_2(A)$. Then the image of $R'$ in $\Aut(H_A/2)$ is
the product of $g$ copies of $SL_2(\F_2)$, and as an $R'$-module, $H_A$ is the
direct sum of $g$-copies of the 2-dimensional vector space over $\F_2$, from
which it follows that $H^0(R_A,H_A/2)$ vanishes, as claimed. Thus
$H^1(R_A,H_A)$ vanishes.

The case $V = L_A/H_A$ follows from the previous case as $L_A/H_A$ can be
imbedded in $L_A$ as an $R_A$-module. The mapping is defined by
$$
x \wedge y \wedge z + H_A \mapsto
(g-1) x \wedge y \wedge z - q \wedge c(x \wedge y \wedge z)
$$
where $q$ denotes the symplectic element and $c$ is the contraction defined
above. It is not difficult to show that this mapping is well defined and that
its composition
$$
L_A/H_A \to L_A \to L_A/H_A
$$
with the quotient mapping is $g-1$ times the identity. Since $L_A/H_A$ is
torsion free and the sequence
$$
0 \to L_A/H_A \to L_A \stackrel{c}{\to} H_A \to 0
$$
is exact, it follows (by tensoring with $\F_2$) that $(L_A/H_A)/2$ injects
into $L_A/2$ from which it follows that $H^0(R_A,L_A/H_A)$ vanishes, as 
required.
\end{proof}

\section{Topological Computations}

\subsection{Johnson's Theorem}
\label{sec:john_thm}
Suppose that $(S,x)$ is a pointed, compact oriented surface of genus $g \ge 2$.
As in previous sections, we denote $H_1(S,A)$  by $H_A$ where $A$ is a ring
such as $\Z,\Zl,\Q,\Ql$, and $(\Lambda^3 H_A)(-1)$ by $L_A$. 

The Torelli group $T_{S,x}$ is, by definition, the kernel of the natural
homomorphism $\G_{S,x} \to \Sp(H_A)$. By the preceding section, one has
a Johnson homomorphism
$$
\tau_{S,x} : H_1(T_{S,x}) \to L_\Z.
$$

\begin{theorem}[\cite{johnson:h1}]
\label{thm:johnson}
If $g\ge 3$, the Johnson homomorphism $\tau_{S,x}$ is surjective and has finite
kernel of exponent $2$. \qed
\end{theorem}

\subsection{Consequences}
Suppose that $V$ is an $\Sp(H_\Z)$-module. The spectral sequence
$$
H^s(\Sp(H_\Z),H^t(T_{S,x})\otimes V) \implies H^{s+t}(\G_{S,x},V)
$$
associated to the group extension
$$
1 \to T_{S,x} \to \G_{S,x} \to \Sp(H_\Z) \to 1
$$
yields the exact sequence
\begin{multline}
\label{seqce}
0 \to H^1(\Sp(H_\Z),V) \to H^1(\G_{S,x},V) \cr
\to H^0(\Sp(H_\Z),\Hom_\Z(H_1(T_{S,x}),V)) \to H^2(\Sp(H_\Z),V).
\end{multline}

The modules $H_\Q$ and $L_\Q$ are rational representations of the algebraic
group $\Sp(H_\Q)$. The representations $H_\Q$ and $L_\Q/H_\Q$ are irreducible
and are the first and third fundamental representations, respectively.

\begin{proposition}[{\cite[Prop.~5.2]{hain:normal}}]
\label{prop:cohomology-computation}
If $g\ge 3$ and $V$ is a finite dimensional irreducible representation of
$\Sp(H_\Q)$, then $H^1(\G_{S,x},V)$ vanishes unless $V$ is isomorphic to either
$H_\Q$ or $L_\Q/H_\Q$. In these cases we have
$$
H^1(\G_{S,x},H_\Q) \cong \Q \text{ and } H^1(\G_{S,x},L_\Q/H_\Q) \cong \Q.
$$
\end{proposition}

\begin{proof}
This follows immediately from Johnson's Theorem (Thm.~\ref{thm:johnson}), the
decomposition (\ref{eq:decomp}), the sequence (\ref{seqce}), Schur's Lemma, and
the fact that $H^j(\Sp(H_\Q),V)$ vanishes for $j=1,2$ by a Theorem of Borel
\cite{borel:twisted} (cf.\ \cite[Thm.~3.2]{hain:torelli}).
\end{proof}

The following computations are equivalent to those of Morita in \cite{morita}.

\begin{corollary}
\label{cor:integral_comp}
If $g\ge 3$, we have
$H^1(\G_{S,x},H_\Z) \cong \Z, \quad H^1(\G_{S,x},L_\Z/H_\Z) \cong \Z$
\text{ and } $H^1(\G_{S,x},L_\Z) \cong \Z^2$.
\end{corollary}

\begin{proof}
It follows from Proposition~\ref{prop:vanish}, using the exact sequence
(\ref{seqce}), that the mapping
$$
H^1(\G_{S,x},V) \to H^0(\Sp(H_\Z),\Hom_\Z(H_1(T_{S,x}),V))
$$
is injective when $V = H_\Z$, $L_\Z$ and $L_\Z/H_\Z$. Since the right hand
group is torsion free, it follows that for these $V$, $H^1(\G_{S,x},V)$ is a
torsion free $\Z$-module of rank equal to the dimension of
$H^1(\G_{S,x},V\otimes\Q)$. The result now follows from the previous result.
\end{proof}

Denote by $\mhat^\an$ the unique element of $H^1(\G_{S,x},L_\Z)$ whose image in
$$
\Hom_{\Sp(H_\Z)}(H_1(T_{S,x}),L_\Z)
$$
is twice the Johnson homomorphism $2\tau_{S,x}$. (Cf.\
Proposition~\ref{prop:univ_class}.)

\subsection{A monodromy computation}
\label{sec:monod_comps}

In this paragraph, we give a detailed proof of the following assertion which is
crucial in this paper, and whose proof is only sketched in \cite{hain:normal}.

\begin{proposition}
\label{prop:equality}
The classes $\mhat^\an$ and $\mu^\an(\cC_g^1,\xihat)$ of
$H^1(\M_g^{1\,\an},\L_\Z)$ are equal.
\end{proposition}

\begin{proof}
Denote $\mu^\an(\cC_g^1,\xihat)$ by $\mu^\an$. It follows from the injectivity
of $H^1(\G_{S,x},L_\Z) \to H^0(S,\Hom_\Z(H_1(T_{S,x}),L_\Z))$ that to prove the
result, we need only check the equality of $\mhat^\an$ and $\mu^\an$ in
$\Hom(H_1(T_{\Canx}),L_\Z)$. We do this by showing that the image of $\mu^\an$
is $2\tau_{S,x}$.

Fix a pointed compact Riemann surface $(S,x)$ as a base point of
$\M_g^{1\,\an}$. Each $\phi \in T_{S,x}$ can be represented by an orbifold
imbedding
$$
\Phi : (S^1,1) \to \big(\M_g^{1\,\an},[C,x]\big).
$$
The pullback of the universal curve along $\Phi$ is the mapping torus
$$
M(\phi) = (S \times [0,1])/\{(\phi(s),0) \sim (s,1): s \in S\}
$$
which projects to $S^1$, which we identify with the unit interval with its ends
identified. The universal section $\xihat$ pulls back to the section that takes
$t\in [0,1]/\{0\sim 1\}$ to $(x,t) \in M(\phi)$.

As a real manifold, the jacobian of $S$ is $H_1(S,\R/\Z)$, which we shall
denote by $J$. Since $\phi$ acts trivially on $H_1(S)$, the bundle of jacobians
has a natural trivialization over $S^1$. The imbedding of the universal curve
into the universal jacobian given by $\xihat$ restricts to an
imbedding\footnote{Constructed explicitly in \cite{hain:normal}.}  $M(\phi)
\to S^1\times J$ over $S^1$. Identify $M(\phi)$ with its image. The cycle
$\cC_{g,\xihat} - \cC_{g,\xihat}^-$ restricts to the topological 3-cycle
$$
Z := M(\phi) - M(\phi)^-
$$
in $S^1\times J$ over $S^1$. Denote the fiber of $Z$ over $t \in S^1$ by $Z_t$.

By duality, the local system used to define the invariant
$\mu^\an(\cC_g^1,\xihat)$ pulls back to the local system over $S^1$ whose fiber
over $t$ is $H_3(J,Z_t)$. This is an extension
$$
0 \to H_3(J) \to H_3(J,Z_t) \to \Z \to 0
$$
in which a relative cycle whose boundary is $Z_t$ projects to 1. The monodromy
of this local system is trivial on $L_\Z$ and $\Z$, and is therefore given by a
homomorphism $\Z \to H_3(J)$, which is, in turn, given by an element of
$H_3(J)$. We will show that this invariant is
$$
2p_\ast[M(\phi)] \in H_3(J) \cong L_\Z
$$
where $p : S^1 \times J \to J$ is the projection. The result will then follow
from Johnson's assertion \cite{johnson:survey} that
$$
p_\ast[M(\phi)] = \tau_{S,x}(\phi),
$$
which is proved in \cite{hain:normal}.

We first work in the pullback $[0,1]\times J \to [0,1]$ along $[0,1] \to S^1$.
Choose a 3-chain $W_0$ in $J$ such that $\partial W_0 = Z_0$. Let $Z_{[0,t]}$
be the restriction of $Z$ to $[0,t]$. Note that
$$
\partial Z_{[0,t]} = Z_t - Z_0.
$$
If we set $\Wtilde_t = W_0 + Z_{[0,t]}$, then $\partial \Wtilde_t = Z_t$. To
get a 3-chain in $J$, set $W_t = p_\ast \Wtilde_t$, where $p : [0,1]\times J
\to J$ is the projection. Then $W_t$ is a 3-chain in $J$ such that $\partial
W_t = Z_t$. The monodromy is now the class of
$$
W_1 - W_0 = p_\ast(\Wtilde_1 - \Wtilde_0) = p_\ast(Z_{[0,1]})
= p_\ast(M(\phi)) - p_\ast(M(\phi)^-)
$$
in $H_3(J)$, which is $2p_\ast([M(\phi)])$, as claimed.
\end{proof}

\subsection{Generators}

In this section we show that several of the universal classes defined in
Section~\ref{sec:univ_classes} turn out to be integral generators of the
cohomology groups in which they lie. The following result is closely related to
computations of Morita in \cite{morita} and \cite{morita:jac} where he
identifies generators of $H^1(\G_{S,x},V)$ where $V$ is $H_\Z$ and $L_\Z/H_\Z$.

\begin{proposition}
\label{prop:cohomology-generator}
If $g \ge 3$, then
\begin{enumerate}
\item $\nutilde^\an(\cC_g^1,\xihat)$ freely generates
$H^1(\M_g^{1\,\an},\L_\Z/\H_\Z)$;
\item $\k^\an(\cC_g^1,\xihat)$ freely generates $H^1(\M_g^{1\,\an},\H_\Z)$;
\item $\nutilde(\cC_g^1,\xihat)$ freely generates $\Het^1(\M_g^1\otimes
\Kbar,\L_\Zl/\H_\Zl)$;
\item $\k(\cC_g^1,\xihat)$ freely generates
$\Het^1(\M_g^1\otimes\Kbar,\H_\Zl)$.
\end{enumerate}
Moreover, $\nutilde^\an(\cC_g^1,\xihat)$ and $\nutilde(\cC_g^1,\xihat)$
correspond under the comparison isomorphism
$$
H^1(\M_g^{1\,\an},\L_\Z/\H_\Z)\otimes\Zl
\cong \Het^1(\M_g^1\otimes \Kbar,\L_\Zl/\H_\Zl)
$$
and $\k^\an(\cC_g^1,\xihat)$ and $\k(\cC_g^1,\xihat)$ correspond under the
comparison isomorphism
$$
H^1(\M_g^{1\,\an},\H_\Z)\otimes\Zl \cong \Het^1(\M_g^1\otimes\Kbar,\H_\Zl)
$$
\end{proposition}

\begin{proof}
It suffices to prove the analytical version, as the \'etale version then
follows via comparison theorems. We therefore restrict ourselves to the
analytical case. By Corollary~\ref{cor:integral_comp}, both groups are torsion
free groups of rank one. So it suffices to show that these classes are both
non-zero and are not divisible by a positive integer.

By Proposition~\ref{prop:equality}, to prove (i), it suffices to show that the
image of $\mhat^\an$ in $H^1(\G_{S,x},L_\Z/H_\Z)$ is not divisible. But this
is true because the image of $\mhat^\an$ under the inclusion
$$
H^1(\G_{S,x},L_\Z/H_\Z) \hookrightarrow \Hom_{\Sp(H_\Z)}(H_1(T_{S,x}),L_\Z/H_\Z)
$$
is $2\tau_{S,x}$, and because Morita's computations \cite[p.~220]{morita}
show that $\tau_{S,x}$ is not in the image of this map.

We now prove that $\kappa^\an(\cC_g^1,\xihat)$ is not divisible. For
convenience, we denote it by $\kappa^\an$. We have surjections
$$
H_1(T_{S,x}) \to L_\Z \to H_\Z.
$$
By \cite[Cor.~6.7]{pulte}, the image of $\mu^\an$ under the induced mapping
$H^1(\M_g^{1\,\an},\L_\Z) \to H^1(\M_g^{1\,\an},\H_\Z)$ is $\kappa^\an$.  Since
the above mapping is surjective, and since $\mu^\an$ goes to $2\tau_{S,x}$, it
follows that if $\kappa^\an/d$ is integral and $d > 1$, then $d=2$. But this is
impossible as it would imply that the section $(2g-2)x - K_C$ of the relative
jacobian $\J \to \M_g^{1\,\an}$ is divisible by 2, which would imply that there
is a square root of $K_C$ defined over $\M_g^1$. If $-2\sigma = (2g-2)x - K_C$,
then
$$
K_C = 2(\sigma + (g-1)x).
$$
That is, there is a theta characteristic defined over $\M_g^{1\,\an}$. But
there is no such theta characteristic. It follows that $\kappa^\an$ generates
$H^1(\M_g^{1\,\an},\H_Z)$.

The statement that $\kappa^\an$ is not divisible also follows from Pulte's
result and \cite[Prop.~6.4]{morita:jac}.
\end{proof}

The wedging with the symplectic element map induces a homomorphism
$$
q_\ast : H^1(\M_g^{1\,\an},\H_\Z) \to H^1(\M_g^{1\,\an},\L_\Z).
$$
Likewise in the \'etale case:

\begin{corollary}
If $g\ge 3$, then
$$
H^1(\M_g^{1\,\an},\L_\Z) =
\Z\, q_\ast\kappa^\an(\cC_g^1,\xihat) \oplus \Z\, \mu^\an(\cC_g^1,\xihat).
$$
\end{corollary}

\begin{proof}
This follows immediately from the previous result using the exact sequence of
cohomology groups
$$
\cdots \to H^1(\M_g^{1\,\an},\H_\Z) \to H^1(\M_g^{1\,\an},\L_\Z)
\to H^1(\M_g^{1\,\an},\L_\Z/\H_\Z) \to \cdots
$$
associated to the short exact sequence
$0 \to \H_\Z \stackrel{q}{\to} \L_\Z \to \L_\Z/\H_\Z \to 0$
of coefficients.
\end{proof}

\begin{corollary}
\label{cor:isom}
If $g\ge 3$, the natural homomorphism 
$$
\Het^1(\M_g^1,\L_\Zl) \to \Het^1(\M_g^1\otimes \Kbar,\L_\Zl)
$$
is an isomorphism. Moreover, both groups are generated freely by
$\mu(\cC_g^1,\xihat)$ and $q_\ast\k(\cC_g^1,\xihat)$.
\end{corollary}

\begin{proof}
The fact $\Het^1(\M_g^1\otimes \Kbar,\L_\Zl)$ is freely generated by
$\mu(\cC_g^1,\xihat)$ and $q_\ast\k(\cC_g^1,\xihat)$ follows from the previous
result via comparison theorems. Since these classes are restrictions of classes
in $\Het^1(\M_g^1,\L_\Zl)$, they are Galois invariant. Using the spectral
sequence\footnote{Note that this spectral sequence uses the definition of
\'etale cohomology given in \cite{jannsen}. Although the definitions of
$\Het^\dot$ in \cite{jannsen} and \cite{milne} differ in general, they coincide
for $\Het^1$.}
$$
H^s(G_K,\Het^t(\M_g^1\otimes \Kbar,\L_\Zl)) \implies \Het^{s+t}(\M_g^1,\L_\Zl).
$$
constructed in \cite{jannsen} and the fact that $H^0(\M_g^1\otimes
\Kbar,\L_\Zl)$ vanishes, which follows directly from
Proposition~\ref{prop:vanish}, it follows that
\begin{multline*}
\Het^1(\M_g^1,\L_\Zl) \cong H^0(G_K,\Het^1(\M_g^1\otimes \Kbar,\L_\Zl)) \cr
\cong \Het^1(\M_g^1\otimes \Kbar,\L_\Zl) \cong \Zl\,q_\ast\k(\cC_g^1,\xihat)
\oplus\Zl\,\mu(\cC_g^1,\xihat).
\end{multline*}
\end{proof}

\begin{remark}
\label{rem:projection}
It is worth noting that the image of $\mu^\an(\cC_g^1,\xihat)$ under
the homomorphism
$$
H^1(\M_g^{1\,\an},\L_\Z) \to H^1(\M_g^{1\,\an},\H_\Z)
$$
induced by the contraction $c:\L_\Z \to \H_\Z$, defined in
Section~\ref{sec:torsion-comps}, is $\k^\an(\cC_g^1,\xihat)$. A proof can be
found in \cite[Cor.~6.7]{pulte}. It thus follows from comparison theorems that
the image of $\mu(\cC_g^1,\xihat)$ under the map
$$
\Het^1(\M_g^1,\L_\Zl) \to \Het^1(\M_g^1,\H_\Zl)
$$
induced by the contraction $c$ is $\k(\cC_g^1,\xihat)$.
\end{remark}

For cohomology with $\Ql$-module coefficients, the results of this section may
be summarized by the following easily proved extension of
Proposition~\ref{prop:cohomology-computation}.

\begin{theorem}
\label{thm:coho_comp}
If $g\ge 3$ and $V$ is a finite dimensional irreducible representation of
$\GSp(H_\Ql)$, then $\Het^1(\M_g^1,\V)$ vanishes unless $V$ is isomorphic to
$H_\Ql$, $L_\Ql/H_\Ql$ or some $\Ql(n)$. In these cases we have
$$
\Het^1(\M_g^1,\H_\Ql) \cong \Ql\k(\cC_g^1,\xihat) \text{ and }
\Het^1(\M_g^1,\L_\Ql/\H_\Ql) \cong \Ql\nu(\cC_g^1)
$$
and $\Het^1(\M_g^1,\Ql(n)) \cong H^1(G_K,\Ql(n))$. \qed
\end{theorem}

\section{Proof of the Density Theorems}
\label{sec:density}

We begin by relating the notation of the Introduction to that of the preceding
sections. Throughout the remainder of this paper, $K$ is a subfield of $\C$ and
$(C,\xi)$ is a pointed smooth projective curve over $K$ of genus $g \ge 3$. As
in the introduction, we denote the $\Ql$-Lie algebra of $\pi_1^\un(\Canx)$ by
$\p(\Cx)$. Since the exponential mapping $\p(\Cx) \to \pi_1^\un(\Cx)$ is an
isomorphism of provarieties, and since, this is an isomorphism of prounipotent
groups if we endow $\p(\Cx)$ with the product given by the
Baker-Campbell-Hausdorff formula, the exponential mapping induces a natural
group isomorphism
$$
\Aut \p(\Cx) \stackrel{\simeq}{\longrightarrow} \Aut \pi_1^\un(\Canx).
$$
The homomorphisms $\rho_{\Cx}$ and $\theta_\Cx$ defined in the Introduction
fit into the commutative diagram
$$
\xymatrix{
& \Aut \pi_1^\un(\Canx) \ar[dd]^\cong
&\cr G_K \ar[ur]^{\rho^\un_{\Cx}} \ar[rd]_{\rho_\Cx} & & \G_\Canx
\ar[lu]_{\theta^\un_{\Cx}} \ar[ld]^{\theta_{\Cx}}\cr & \Aut \p(\Cx) &
}
$$
In short, we have removed the decoration ``$\un$'' when the target is $\Aut
\p(\Cx)$ after making the natural identifications.

Since $K$ is embedded in $\C$, we can naturally regard $[\Cx]$ as a
$\C$-rational point of $\M_g^1$ to which corresponds the homomorphism
$$
\rhohat_\Cx : \pi_1^\alg(\M_g^1,[\Cx]) \to \Aut\p(\Cx)
$$
that extends both $\theta_\Cx$ and $\rho_\Cx$:
\begin{equation}
\label{diag:compat}
\xymatrix{
G_K \ar[rrd]_{\rho_\Cx} \ar[rr]^(0.4){[\Cx]_\ast} & &
\pi_1^\alg(\M_g^1,[\Cx])\ar[d]^{\rhohat_\Cx}  & &
\G_\Canx\ar[lld]^{\theta_\Cx}\ar[ll] \cr
& & \Aut \p(\Cx)
}
\end{equation}

\subsection{Strategy}
\label{sec:strategy}
We shall prove each of the 3 statements of Theorem~\ref{th:main-aut}
together with the following:
\begin{quote}
(iv)
{\it The Zariski closure of the image of $\rho_{\Cx}$ in $\Aut\p(\Cx)$  is
identical with that of the image of $\rhohat_\Cx$.}
\end{quote}
We shall do this by proving that (i) $\Rightarrow$ (ii) $\Rightarrow$ (iii)
$\Rightarrow$ (iv) $\Rightarrow$ (i). The implications (iv) $\Rightarrow$ (i)
$\Rightarrow$ (ii)  are obvious.  In this section, we shall prove  (ii)
$\Rightarrow$ (iii) and (iii) $\Rightarrow$ (iv) to complete the proof.

\subsection{Proof of (ii) $\Rightarrow$ (iii)}
\label{2=>3}
First, if the closure of $\im\rho_\Cx^2$ contains $\im\theta_\Cx^2$, then the
closure of the image of $\rho_\Cx^1$ contains $\im\theta_\Cx^1$. The image of
$\theta_\Cx^1$ is $\Sp(H_\Z)$, and therefore has Zariski closure $\Sp(H_\Ql)$
in $\GSp(H_\Ql)$. Since the composite
$$
G_K \stackrel{\rho_1}{\to} \GSp(H_\Ql) \stackrel{\chi}{\to} \Gm
$$
is the $\l$-adic cyclotomic character, and since it is assumed to have infinite
image, it follows that if the Zariski closure of $\im\rho_\Cx^1$ contains
$\im\theta_\Cx^1$, then its Zariski closure is all of $\GSp(H_\Ql)$.

It is well-known and easy to show that the Lie algebra $\p(\Cx)/L^3$ is
(unnaturally) isomorphic to the graded Lie algebra
$$
H_\Ql \oplus \Lambda^2 H_\Ql/q
$$
where $q$ spans the $\Sp(\H_\Ql)$ invariants in $\Lambda^2 H_\Ql$. The
bracket is given by the mapping
$$
H_\Ql\otimes H_\Ql \to \Lambda^2 H_\Ql/q.
$$
that takes $a\otimes b$ to $a\wedge b \mod q$. Using this description, it
is easy to see that $\Aut \big[\p(\Cx)/L^3\big]$ is an extension
$$
0 \to \Hom(H_\Ql,\Lambda^2 H_\Ql/q) \to \Aut \big[\p(\Cx)/L^3\big]
\to \GSp(H_\Ql) \to 1.
$$
The kernel is abelian. When $g\ge 3$, the kernel decomposes as an $\GSp$-module
into three irreducible components:
$$
\Hom(H_\Ql,\Lambda^2 H_\Ql/q) \cong H_\Ql \oplus L_\Ql/H_\Ql \oplus U(-1)
$$
where $U$ is the irreducible representation corresponding to the partition
$[2,1]$.\footnote{Each partition $\lambda$ of a non-negative integer $m$ into
$\le g$ pieces corresponds to an irreducible representation of $\Sp(H)$. Each
of these can be lifted uniquely to a representation $V_\lambda$ of $\GSp(H)$
where the scalar matrix $aI$ acts as $a^m\id_V$. Every other representation of
$\GSp(H)$ can be obtained from such a $V_\lambda$ by tensoring with the one
dimensional representation $\Z(n)$ for some $n$. We shall denote the weight of
$V_\lambda(n) := V_\lambda\otimes \Z(n)$ by $\lambda(n)$. This has weight
$-m-2n$, which means that the scalar matrix $aI$ acts as $a^{m+2n}\id_V$ on
it.} Each component has weight $-1$ as a $\GSp(H_\Ql)$-module.

To proceed, we need to recall an elementary fact about group cohomology.
Suppose that $R$ is a reductive affine algebraic group over $\Ql$ and that $V$
is an irreducible $R$-module. Suppose $\G$ is a profinite group and that $r:\G
\to R(\Ql)$ is a continuous Zariski dense representation. Each lift of $r$ to a
homomorphism $\rtilde : \G \to R\ltimes V$, determines a cohomology class
$$
c \in H^1(\G,V)
$$
by taking the class of the cocycle $f$ obtained by taking the second component
of $\rtilde$:
$$
\rtilde(g) = (r(g),f(g)) \in R\ltimes V.
$$
The class $c$ determines $\rtilde$ up to conjugation by an element of $V$. The
following result is easily proved. Details are left to the reader.

\begin{lemma}\label{lem:cohom-cond}
The Zariski closure of the image of the lift $\rtilde$ contains $V$ if and only
if $c$ is not zero. \qed
\end{lemma}

Set
$$
E = \Hom(H_\Ql,\left[\Lambda^2 H_\Ql\right]/q).
$$
This decomposes as the sum of one copy each of $V_\lambda$ for $\lambda \in
\{[1],\, [1^3](-1),\, [2,1](-1)\}$. The group $\Aut \big[\p(\Cx)/L^3\big]$ is
an extension
$$
0 \to E \to \Aut \big[\p(\Cx)/L^3\big] \to \GSp(H_\Ql) \to 1.
$$
Since the diagram
$$
\xymatrix{
G_K \ar[rrd]_{\rho_\Cx^2} \ar[rr]^(0.4){[\Cx]_\ast} & &
\pi_1^\alg(\M_g^1,[\Cx])\ar[d]^{\rhohat_\Cx^2}  & &
\G_\Canx\ar[lld]^{\theta_\Cx^2}\ar[ll] \cr
& & \Aut \big[\p(\Cx)/L^3\big]
}
$$
commutes, the construction above gives characteristic classes
$$
c_\lambda \in H^1(G_K,E),\ \ctilde_\lambda \in H^1(\pi_1^\alg(\M_g^1,[\Cx]),E),\
\ctilde^\an_\lambda \in H^1(\G_\Canx,E)
$$
for each $\lambda \in \{[1],\, [1^3](-1),\, [2,1](-1)\}$. Note that
$\ctilde_\lambda$ maps to $c_\lambda$ and $\ctilde^\an_\lambda$ under the
induced maps
$$
H^1(G_K,E) \leftarrow H^1(\pi_1^\alg(\M_g^1,[\Cx]),E)
\rightarrow H^1(\G_\Canx,E).
$$
It follows from \cite[Cor.~9.2]{hain:torelli} that the Zariski closure of the
image of $\theta_\Cx^2$ contains $V_{[1]}$ and $V_{[1^3](-1)}$, but does not
intersect $V_{[2,1](-1)}$. It follows that the Zariski closure of $\rho_\Cx^2$
contains $\im\theta_\Cx^2$ if and only if $c_{[1]}$ and $c_{[1^3](-1)}$ are
both non-zero. To complete the proof that Proof of (ii) $\Rightarrow$ (iii) we
will show that $c_{[1]}$ is a non-zero multiple of $\kappa(\Cx)$ and
$c_{[1^3](-1)}$ is a non-zero multiple of $\nu(C)$.

To prove this, note that \cite[Cor.~9.2]{hain:torelli} implies that
$\ctilde^\an_{[1]}$ and $\ctilde^\an_{[1^3](-1)}$ are both non-zero and
therefore, by  Proposition~\ref{prop:cohomology-generator}, $\ctilde^\an_{[1]}$
is a non-zero multiple of $\kappa^\an(\cC_g^1,\xihat)$ and
$\ctilde^\an_{[1^3](-1)}$ is a non-zero multiple of
$\nutilde^\an(\cC_g^1,\xihat)$. Another application of
Proposition~\ref{prop:cohomology-generator} implies that $\ctilde_{[1]}$ is a
non-zero multiple of $\kappa(\cC_g^1,\xihat)$ and that $\ctilde_{[1^3](-1)}$ is
a non-zero multiple of $\nutilde(\cC_g^1,\xihat)$.
Propositions~\ref{prop:K-specialization} and \ref{prop:N-specialization} now
imply that $c_{[1]}$ is a non-zero multiple of $\k(\Cx)$ and that
$c_{[1^3](-1)}$ is a non-zero multiple of $\nu(C) = \nutilde(\Cx)$, as claimed.

\subsection{Weighted completion}
\label{sec:completion}

To prove (iii) $\Rightarrow$ (iv) we will need the tool of weighted completion.

We begin with a brief review of the theory of weighted completion of a
profinite group. More details (including complete proofs) can be found in
\cite{hain-matsumoto} and \cite{hain-matsumoto:exp}.

Suppose that $F$ is a topological field, $R$ a reductive affine algebraic group
defined over $F$, that $\G$ is a profinite group, and that $r:\G \to R(F)$ is a
continuous homomorphism with Zariski dense image. In addition, we fix a central
cocharacter $w:\Gm \to R$, that is, a homomorphism from $\Gm$ whose image lies
in the center of $R$.

A $\G$-module is a finite dimensional $F$-vector space on which $\G$ acts
continuously. Let $F(m)$ denote a one-dimensional $F$-vector space on which
$\Gm$ acts by $m$-th power multiplication. An $R$-module is of pure weight $m$
if it is isomorphic to a sum of copies of $F(m)$ as a $\Gm$-module via $w$.
Schur's lemma implies that each irreducible $R$-module is  pure of weight $m$
for some integer $m$.

A {\it weighted $\G$-module} is a finite dimensional $\G$-module $M$ with an
increasing filtration (necessarily unique)
$$
0 = W_nM \subseteq W_{n+1}M \subseteq \cdots \subseteq W_{N-1}M
\subseteq W_N M = M, 
$$
by $\G$-submodules, where, for each integer $m$, the weight graded quotient
$\Gr^W_m M := W_m M / W_{m-1}M$ is an $R$-module of weight $m$ in such a way
that the action of $\G$ on it factors through this $R$ action:
$$
\G \stackrel{r}{\to} R(F) \to \Aut \Gr^W_m V.
$$
One fact that makes weighted $\G$-modules useful is:

\begin{proposition}[\cite{hain-matsumoto}]
The functor $\Gr^W_\dot$ is exact on the category of weighted $\G$-modules.
\qed
\end{proposition}

The category of weighted $\G$-modules is tannakian with fiber functor the
functor that takes a weighted $\G$-module to its underlying vector space.

The {\it weighted completion} of $\G$ with respect to $r$ and $\w$ is the
fundamental group of the tannakian category of weighted $\G$-modules with
respect to the fiber functor above. We shall denote it by $\cG(\G,r)$. It is a
pro-algebraic group over $F$, which is an extension
$$
1 \to \U(\G,r) \to \cG(\G,r) \to R \to 1
$$
of $R$ by a prounipotent group $\U(\G,r)$. There is a natural homomorphism
$\rtilde : \G \to \cG(\G,r)(F)$ that lifts $r$, and is easily seen to be
Zariski dense. The action of $\G$ on a weighted $\G$-module factors through
this homomorphism.

\subsubsection*{The setting}
Fix a pointed smooth projective curve $(C,\xi)$ of genus $\ge 3$ over the
subfield $K$ of $\C$. Denote the dual of $\Het^1(C\otimes \Kbar,\Ql)$ by
$H_\Ql$. In the rest of this section, $F$ will be $\Ql$ and $R$ will be
$\GSp(H_\Ql)$. The central cocharacter is the one that takes $a \in \Gm$ to
$a^{-1}\id$. The group $\G$ is any that acts naturally on $C$, such as $G_K$ or
$\pi_1^\alg(\M_g^1,[C_x,x])$, where $x$ is a geometric point of $C$ lying over
$\xi$. We suppose that the induced representation $r : \G \to \GSp(H_\Ql)$ is
Zariski dense. We thus have the weighed completion $\cG(\G,r)$.

\begin{proposition}
\label{prop:weight-act}
In this setting, $\p(\Cx)$ is the inverse limit of weighted $\G$-modules. The
action $\G \to \Aut \p(\Cx)$ thus induces a homomorphism $\cG(\G,r) \to
\Aut \p(\Cx)$.
\end{proposition}

\begin{proof}
Define a weight filtration on $\p(\Cx)$ by
by
$$
W_{-m} \p(\Cx) = L^m\p(\Cx).
$$
Since the lower central series of $\p(\Cx)$ is characteristic, it is  preserved
by $\G$. Since $H_\Ql$ has weight $-1$, and since the iterated bracket mapping
$$
H_\Ql^{\otimes m} \to \Gr^W_{-m}\p(\Cx),\qquad x_1\otimes \cdots \otimes x_m
\mapsto [x_1[x_2[\cdots [x_{m-1},x_m]]\cdots]]
$$
is surjective and $\G$-invariant, it follows that $\Gr^W_{-m}\p(\Cx)$ is a
$\GSp(H_\Ql)$-module of weight $-m$ and that each $\p(\Cx)$ is the inverse
limit of the weighted $\G$-modules $\p(\Cx)/W_{-m}$.
\end{proof}

We conclude this paragraph with a surjectivity criterion which will be used
in the proof of (iii) $\Rightarrow$ (iv).

Let $F$ be a topological field and $\G_j$ ($j=1,2$) be two topological groups.
Suppose that $R$ is a reductive group over $F$ with central cocharacter $w :
\Gm \to R$ as above. Suppose also that we have two Zariski dense
representations $r_j:\G_j \to R(F)$.  We then have the two corresponding
weighted completions by
$$
\rhotilde_j : \G_j \to \cG(\G_j) \quad j=1,2.
$$

A homomorphism $\phi : \G_1 \to \G_2$ that commutes with the projections to
$R(F)$ induces a homomorphism $\phi_\cG : \cG_1 \to \cG_2$ that commutes with
the projections to $R$.

\begin{lemma}
\label{lem:surj}
The homomorphism $\phi_\cG : \cG(\G_1) \to \cG(\G_2)$ is  surjective if and
only if the induced mappings
$$
H^1(\G_2, V_\alpha) \to H^1(\G_1, V_\alpha)
$$
are injective for all negatively weighted irreducible representation $V_\alpha$
of $R$.
\end{lemma}

\begin{proof}
Denote the prounipotent radical of $\cG(\G_j)$ by $\U_j$, $j=1,2$. Since the
two $\cG(\G_j)$ have the same reductive quotient, $\phi_\cG$ is surjective if
and only if its restriction
$$
\phi_\U : \U_1 \to \U_2
$$
is surjective. This occurs if and only if the induced mapping $d\phi_\U : \u_1
\to \u_2$ on Lie algebras is surjective. Since the $\u_j$ are pronilpotent,
this occurs if and only if the induced mapping
$$
\phi_{\U,\ast} : H_1(\u_1) \to H_1(\u_2)
$$
on abelianizations is surjective. The $H_1(\u_j)$ are weighted modules and
$\phi_{\U,\ast}$ preserves the weight filtration (cf. \cite[\S\S 3.3,
3.4]{hain-matsumoto}). By the exactness of $\Gr^W_\dot$, $\phi_{\U,\ast}$ is
surjective, if and only if
$$
\Gr^W_\dot \phi_{\U,\ast} : \Gr^W_\dot H_1(\u_1) \to \Gr^W_\dot H_1(\u_2)
$$
is surjective. Dualizing, this is true if and only if
\begin{equation}
\label{induced}
\phi_\U^\ast : H^1(\Gr^W_\dot \u_2) \to  H^1(\Gr^W_\dot \u_1)
\end{equation}
is injective. By \cite[Thm.~4.8]{hain-matsumoto}, there are natural
isomorphisms
$$
H^1(\Gr^W_\dot \u_j) \cong 
\bigoplus_{\alpha} H^1(\G_j,V_\alpha)\otimes V_\alpha^\ast
$$
where $(V_\alpha)_\alpha$ is a set of representatives of the isomorphism
classes of negatively weighted irreducible representations of $R$. The result
follows as the
mapping (\ref{induced}) is induced by the mappings
$$
\phi^\ast : H^1(\G_2,V_\alpha) \to H^1(\G_1,V_\alpha).
$$
\end{proof}

\subsection{Proof of (iii) $\Rightarrow$ (iv)}
Assume the conditions of (iii). Denote the weighted completion of $G_K$ with
respect to
$$
\rho_\Cx^1 : G_K \to \GSp(H_\Ql)
$$
by $\cG(G_K)$, and the weighted completion of $\pi_1^\alg(\M_g^1,[\Cx])$ with
respect to
$$
\rhohat_{\Cx}^1 : \pi_1^\alg(\M_g^1,[\Cx]) \to \GSp(H_\Ql)
$$
by $\cG(\M_g^1)$.

Since the actions
$$
\rho_\Cx : G_K \to \Aut \p(\Cx) \text{ and }
\rho_{\cC_g^1,\xihat} : \pi_1^\alg(\M_g^1,[C_x,\xtilde]) \to \Aut \p(\Cx)
$$
are both negatively weighted, they induce representations
$$
\cG(G_K) \longrightarrow \Aut\p(\Cx) \longleftarrow \cG(\M_g^1).
$$
The homomorphism
$$
[\Cx]_\ast : G_K \to \pi_1^\alg(\M_g^1,[\Cx])
$$
commutes with the projections to $\GSp(H_\Ql)$, and thus induces a homomorphism
$$
(\Cx)_\ast : \cG(G_K) \to \cG(\M_g^1).
$$
The commutativity of the left-hand
diagram below implies the commutativity of the right-hand one:
$$
\xymatrix{
G_K \ar[r]^(0.35){[\Cx]_\ast}\ar[rd] & \pi_1^\alg(\M_g^1,[\Cx]) \ar[d]\cr
& \Aut \p(\Cx)
}
\qquad
\xymatrix{
\cG(G_K) \ar[r]^{(\Cx)_\ast}\ar[rd] & \cG(\M_g^1) \ar[d] \cr
& \Aut \p(\Cx)
}
$$
Since the natural homomorphism $G_K \to \cG(G_K)$ and $\pi_1^\alg(\M_g^1,[\Cx])
\to \cG(\M_g^1)$ have Zariski dense image and, given the commutativity of
(\ref{diag:compat}), (iii) implies (iv) follows directly from the following
more general result.

\begin{theorem}
\label{thm:general_density}
Assuming that the image of $\rho_C^1 : G_K \to \GSp(H_\Ql)$ is Zariski dense,
the homomorphism $(\Cx)_\ast$ is surjective if and only if the classes
$\nu(C)$ and $\kappa(\Cx)$ both have infinite order.
\end{theorem}

\begin{proof}
This follows directly from Proposition~\ref{prop:N-specialization},
Lemma~\ref{lem:surj} and Theorem~\ref{thm:coho_comp}.
\end{proof}

\section{The $\l$-adic Harris-Pulte Theorem}
\label{sec:harris-pulte}

In this section we prove the $\l$-adic Harris-Pulte Theorem. We do this
by proving a more general statement, and deducing it by specialization.

A pointed curve $(\Cx)$ over $K$ and a geometric point $x$ of $\Spec K$
gives rise to a homomorphism
$$
\delta_{C,x} :
\pi_1^\alg(\M_g^1,[C_x,\xtilde]) \to \Aut \pi_1^\prol(C_x,\xtilde),
$$
where $\xtilde = \xi(x)$, and therefore a class
$$
\delta_{C,x}^\ast \mhat^\prol \in H^1(\pi_1^\alg(\M_g^1,[C_x,\xtilde]),L_\Zl)
\cong \Het^1(\M_g^1,\L_\Zl),
$$
where $\mhat^\prol$ is the class defined in \S\ref{sec:char_class}.

This class is easily seen to be independent of the choice of $x$.
We also have the class
$$
\mu(\cC_g^1,\xihat) \in \Het^1(\M_g^1,\L_\Zl)
$$
of the universal cycle $\cC_\xi - \cC_\xi^-$.

The $\l$-adic Harris-Pulte Theorem is a direct consequence of the following
result.

\begin{theorem}
The classes $\delta_{C,x}^\ast\mhat^\prol$ and $\mu(\cC_g^1,\xihat)$ in
$\Het^1(\M_g^1,\L_\Zl)$ are equal.
\end{theorem}

\begin{proof}
First suppose that $K$ is a subfield of $\C$. Denote the algebraic closure of
$K$ in $\C$ by $\Kbar$. Proposition~\ref{prop:equality}, comparison theorems
and the commutativity of the diagram
$$
\begin{CD}
\G_\Canx @>>> \Aut \pi_1(\Canx) \cr
@VVV @VVV \cr
\pi_1^\alg(\M_g^1\otimes\Kbar,[C_x,\xtilde]) @>>{\bar{\delta}_{C,x}}>
\Aut \pi_1^\prol(C_x,\xtilde)
\end{CD}
$$
imply that
$$
\bar{\delta}_{C,x}^\ast \mhat^\prol = \mu(\cC_g^1,\xihat) \in
\Het^1(\M_g^1\otimes\Kbar,\L_\Zl).
$$
The result now follows from Corollary~\ref{cor:isom}. The case of general $K$
follows from the case where $K \subset \C$ by standard techniques of GAGA.
\end{proof}

The Harris-Pulte Theorem now follows by specialization:

\begin{proof}[Proof of Theorem~\ref{th:l-Harris-Pulte}]
As in the statement of the theorem, we suppose that $K$ is a subfield of $\C$
and that $(\Cx)$ is a pointed curve over $K$. Choose a geometric point $x$
of $\Spec K$. We shall identify $\pi_1^\prol(C_x,\xtilde)$ with
$\pi_1^\prol(\Canx)$. The class $m(\Cx)$ is the
pullback of the universal class $\mhat^\prol$ under the natural homomorphism
$$
\rho_{C,x}^\prol : G_K \to \Aut \pi_1^\prol(\Canx).
$$
Since the diagram
$$
\xymatrix{
\Het^1(\M_g^1,\L_\Zl) \ar[r]^{[C,x]^\ast} & H^1(G_K,L_\Zl) \cr
H^1(\Aut \pi_1^\prol(\Canx),L_\Zl) \ar[u]^{\delta_{C,x}^\ast}
\ar[ur]_{\rho_{C,x}^{\prol\ast}}
}
$$
commutes, Proposition~\ref{prop:C-specialization} implies that
$$
m(C,x) =
\rho_{C,x}^{\prol\,\ast} \mhat^\prol = [C,x]^\ast \delta_{C,x}^\ast \mhat^\prol
= [C,x]^\ast \mu(\cC_g^1,\xihat) = \mu(C,x).
$$
\end{proof}

\section{The $\l$-adic Formulation}

As in previous sections, $C$ is a smooth projective curve of genus $\ge 3$
defined over the subfield $K$ of $\C$, $\Kbar$ is its algebraic closure in
$\C$, $\xi$ is a $K$-rational point of $C$, $x$ is the geometric point $\Spec
\Kbar$ of $\Spec K$ and $\xtilde$ is $\xi(x)$.

We have the truncations
$$
\rho_\Cx^{\prol,m} : G_K \to \Aut \big[\pi_1^\prol(\Canx)/L^{m+1}\big]
$$
and
$$
\rhohat_\Cx^{\prol,m} : \pi_1^\alg(\M_g^1,[C_x,\xtilde]) \to
\Aut \big[\pi_1^\prol(\Canx)/L^{m+1}\big]
$$
of the monodromy representations $\rho_\Cx^\prol$ and $\rhohat_\Cx^\prol$
defined in Section~\ref{sec:arith}.

In this section, we prove the following $\l$-adic version of
Theorem~\ref{th:main-aut}.


\renewcommand{\theenumi}{\alph{enumi}}

\begin{theorem}
\label{th:main-aut-4}
If the image of the $\l$-adic cyclotomic character $\chi_\l : G_K \to
\Zl^\times$ is infinite, then the following four conditions are equivalent:
\begin{enumerate}

\item $\im\rho_\Cx^{\prol,m}$ is an open subgroup of $\im\rhohat_\Cx^{\prol,m}$
for all $m\ge 1$.

\item $\im\rho_\Cx^{\prol,2}$ is an open subgroup of
$\im\rhohat_\Cx^{\prol,2}$.

\item The image of $\rho^{\prol,1}_C : G_K \to \GSp(H_\Zl)$ is open, and both
the classes $\k(\Cx)$ in $H^1(G_K,H_\Zl)$  and $\nu(C)$ in
$H^1(G_K,L_\Zl/H_\Zl)$ have infinite order. \item Any of the 3 equivalent
conditions (i), (ii), (iii) of Theorem~\ref{th:main-aut}.
\end{enumerate}
\end{theorem}


\renewcommand{\theenumi}{\roman{enumi}}

Clearly, (a) $\implies$ (b) $\implies$ (c) $\implies$ (iii). To complete the
proof, we now show that (iv) $\implies$ (a), where statement (iv) is given in
Section~\ref{sec:strategy}. This will follow directly from the following three
results. We begin by showing that Zariski density of the action of $G_K$ on
$H_\Ql$ implies that the image of the Galois action on $H_\Zl$ is open.

\begin{lemma}
\label{lem:rational}
If $\G$ is a profinite group and $r: \G \to \GSp(H_\Zl)$ is a continuous
homomorphism whose image in $\GSp(H_\Ql)$ is Zariski dense, then the image of
$r$ is open in $\GSp(H_\Zl)$.
\end{lemma}

\begin{proof}
The image $G$ of $r$ is closed in $\GSp(H_\Zl)$, and therefore an $\l$-adic Lie
group, \cite{serre:LALG}. Its Lie algebra $\g$ is thus an analytic Lie algebra
over $\Ql$. A result of Borel \cite[7.9 Cor.]{borel:GTM} implies that
$[\g,\g]$  is the Lie algebra of an algebraic subgroup $H$ of $\GSp(H_\Ql)$.
This implies, via an argument using the exponential mapping, that  the image of
$[G,G]$ is open in $H(\Ql)$. Since $G$ is Zariski dense in $\GSp(H_\Ql)$,  its
commutator is Zariski dense in $\Sp(H_\Ql)$, and hence  $H$ must contain
$\Sp(H_\Ql)$. But $H$ is contained in the  commutator of $\GSp(H_\Ql)$, and so
must be $\Sp(H_\Ql)$. Since $G$ is dense in $\GSp(H_\Ql)$, the image of the
composite of $r$ with $\GSp(H_\Zl) \to \Zl^\times$ is infinite and therefore
open. The result follows.
\end{proof}

\begin{lemma}
\label{lem:dense-open}
Suppose that $\G$ is a profinite group and that $U$ is a unipotent group over
$\Ql$. If $\G \to U(\Ql)$ is a Zariski dense homomorphism, then the image of
$\G$ is open in $U(\Ql)$.
\end{lemma}

\begin{proof}
See \cite[Lemma~7.5]{hain-matsumoto}.
\end{proof}

The weighted completion $\cG(G_K)$ is an extension
$$
1 \to \U(G_K) \to \cG(G_K) \to \GSp \to 1
$$
of proalgebraic groups, where $\U(G_K)$ is prounipotent.

\begin{lemma}
\label{lem:density-of-kernel}
The image of $\ker\{G_K \to \GSp(\Ql)\}$ in $\U(G_K)$ is Zariski dense.
\end{lemma}

\begin{proof}
By \cite[Corollary~4.5]{hain-matsumoto}, 
it is enough to show the vanishing of 
$$
H^1(\im(G_K), V)
$$
for any negative irreducible $\GSp$-module $V$ over $\Ql$, where $\im(G_K)$
denotes the image of $G_K$ in $\GSp(\Ql)$. By Lemma~\ref{lem:rational},
$\im(G_K)$ contains a non-torsion scalar matrix $aI$. Then, $aI$ is in the
center of $\im(G_K)$, and hence acts trivially on $H^1(\im(G_K),V)$. On the 
other hand, $aI$ acts on $V$ by multiplication by  a negative power of $a$.
This shows that the cohomology  is torsion, hence is trivial over $\Ql$.
\end{proof}

\begin{proof}[Proof of Theorem~\ref{th:main-aut-4}]
Fix a positive integer $m$. Denote by $G$ the image of $\cG(\M_g^1)$ in the
algebraic group $\Aut[\p(\Cx)/L^{m+1}]$. The Zariski density of $\chi_\l$
implies that $G$ is an extension
$$
1 \to U \to G \to \GSp \to 1,
$$
where $U$ is a unipotent group over $\Ql$. Condition (iv) implies that $G$ is
also the image of $\cG(G_K)$ in $\Aut[\p(\Cx)/L^{m+1}]$. What we want to show
is that the image of $G_K$ in $G(\Ql)$ is open.  Denote the kernel of $G_K \to
\GSp(\Ql)$ by $\Delta$.  Since the image of $G_K$ in $\GSp(\Ql)$ is open
(Lemma~\ref{lem:rational}), it suffices to show that $\Delta$ has open image in
$U(\Ql)$. By Lemma~\ref{lem:dense-open}, it suffices  to show that the image of
$\Delta$ is Zariski dense in $U$. But this follows from 
Lemma~\ref{lem:density-of-kernel} and the surjectivity of $\U(G_K) \to U$,
which is a consequence of  the surjectivity of $\cG(G_K) \to G$.
\end{proof}

\section{The Un-pointed Case}
\label{sec:non-pointed}

In this section, we briefly sketch the proof of Theorem~\ref{th:main-out}.
Since the arguments are very similar to those in the pointed case, the details
are left to the reader. We also state its $\l$-adic version and an un-pointed
version of the Harris-Pulte Theorem that does not require the existence of a
rational point. The proofs of these later two results are very similar to those
in the pointed case and are left to the reader.

\subsection{Monodromy representations}
\label{sec:monod}

There are similar constructions to those in Section~\ref{sec:arith} in the
non-pointed case, where $\Aut$ is replaced by $\Out$. Specifically, suppose
that $C \to B$ is a proper smooth family of genus $g$ curves defined over the
field $K$ and that $x$ is a geometric point of $B$. Even though this may not
have a section, there is still a natural homomorphism
$$
\rho_C : \pi_1^\alg(B,x) \to \Out \pi_1^\alg(C_x),
$$
where $C_x$ denotes the geometric fiber of $C$ over $x$. This is constructed
as the composite
$$
\begin{CD}
\pi_1^\alg(B,x) @<{\simeq}<< \pi_1^\alg(C,\tilde{x})/\pi_1^\alg(C_x,\xtilde)
@>{\mathrm{conjugation}}>> \Out \pi_1^\alg(C_x,\xtilde)
\end{CD}
$$
where $\xtilde$ is a geometric point of $C$ lying over $x$. This induces
outer actions
\begin{equation}
\label{eq:monod_nonptd:*}
\rho_C^\ast : \pi_1^\alg(B,x) \to \Out \pi_1^\ast(C_x)
\end{equation}
where $\ast \in \{\alg,\prol,\un\}$.

Denote the moduli stack of smooth projective curves of genus $g>1$ over $K$ by
$\M_g$. It is defined over $\Spec K$. Each pointed curve $C$ defined over an
algebraically closed extension of $K$ determines a geometric point $[C]$ of
$\M_g$; the fiber of the universal curve over this point may be identified with
$C$. Associated to this geometric point is the monodromy representation
\begin{equation}
\label{eq:univ_monod_unptd}
\rhohat_{C} : \pi_1^\alg(\M_g,[C]) \to \Out \pi_1^\alg(C).
\end{equation}
This is called the {\it universal monodromy representation} for the
curve $C$ for reasons that we now recall.

For each proper family $C \to B$ of smooth pointed curves of genus $g$
defined over $K$, there is a unique morphism $[C]: B \to \M_g$ that
classifies the family. It takes the geometric point $x$ of $B$ to $[C_x]$. The
family $C \to B$ is the pull-back of the universal curve $\cC_g \to \M_g$
along $[C]$.

By universality, the fiber of $\cC_g \to \M_g$ over $[C](x)$ can be identified
naturally with $C_x$. Thus, the diagram
$$
\xymatrix{
\pi_1^\alg(B,x)\ar[r]^(0.4){\rho_{C}^\alg}\ar[d]_{[C]_\ast} &
\Out \pi_1^\alg(C_x)\cr
\pi_1^\alg(\M_g^1,[C_x]) \ar[ur]_{\rhohat_{C_x}^\alg}
}
$$
commutes. That is, the monodromy representation of $C \to B$ associated to
the geometric point $x$ of $B$ factors through the universal monodromy
representation (\ref{eq:univ_monod_unptd}) associated to $(C_x,\xtilde)$.

By functoriality, the universal monodromy representation
(\ref{eq:univ_monod_unptd}) also has pro-$\l$ and prounipotent incarnations:
\begin{equation}
\label{univ_monod_nonptd:l}
\rhohat_{C}^\prol : \pi_1^\alg(\M_g,[C]) \to \Out \pi_1^\prol(C)
\end{equation}
and
\begin{equation}
\label{univ_monod_nonptd:un}
\rhohat_{C}^\un : \pi_1^\alg(\M_g,[C]) \to \Out \pi_1^\un(C)
\end{equation}
through which the representations (\ref{eq:monod_nonptd:*}) factor.

Now suppose that $K$ is a subfield of $\C$ and that $C$ is a curve defined over
$K$. Let $\Kbar$ be the algebraic closure of $K$ in $\C$. As in the pointed
case, there is a natural isomorphism
$$
\Ghat_\Can \cong \pi_1^\alg(\M_g\otimes \Kbar,[\Cbar]),
$$
where $\Cbar = C\otimes\Kbar$, which is unique up to conjugation by an
element of $\Aut(\Can)$. Oda's result \cite{oda} implies that $\theta_{\Cbar}$
is the completion of the tautological action
$$
\G_\Can \to \Out \pi_1(\Can).
$$

\subsection{The Johnson class}

Let $\pi$ and $A$ be as in the beginning of Section~\ref{sec:johnson}. Since
the inner automorphisms of $\pi$ act trivially on $H_1(\pi)$, the Johnson
homomorphism induces homomorphisms
$$
\Inn \pi \to H_1(\pi) \to H_1(T) \to L_A.
$$
The image of the composition of these is well known and easily seen to be the
standard copy of $H_A$ in $L_A$. It follows that if we write $\Out \pi$ as an
extension
$$
1 \to OT \to \Out \pi \to \GSp(H_A) \to 1,
$$
where $OT := T/\Inn \pi$, then the Johnson homomorphism $\tau : T \to L_A$
induces a natural, $\GSp(H_A)$-invariant surjection
$$
H_1(OT) \to L_A/H_A.
$$

Arguing as in Section~\ref{sec:char_class}, one can prove the following
analogue of Proposition~\ref{prop:univ_class} for outer actions.

\begin{proposition}
There is a unique class $\nhat \in H^1(\Out\pi, L_A/H_A)$ whose pullback to
$H^1(\Aut \pi, L_A/H_A)$ is the reduction of $\mhat \in H^1(\Aut \pi, L_A)$ mod
$H_A$. \qed
\end{proposition}

In the $\l$-adic case, we shall denote by $\nhat^\prol$ this universal
characteristic class in $H^1(\Out \pi,L_\Zl/H_\Zl)$.

For each curve $C$ over a field $K$ of characteristic 0, define
$$
n(C) \in H^1(G_K,L_\Zl/H_\Zl)
$$
to the the pullback of the universal class $\nhat^\prol$ along the natural
action
$$
G_K \to \Out \pi_1^\prol(C\otimes \Kbar).
$$

In the discrete case, we have Johnson's Theorem, which allows the computation
of various cohomology groups of $\M_g^\an$, and $\M_g$. The Torelli group $T_S$
of a compact oriented surface $S$ of genus $g$ is the kernel of the natural
homomorphism $\G_S \to \Aut H_1(S)$. Johnson's homomorphism for pointed
surfaces induces an $\Sp(H_\Z)$-equivariant homomorphism
$$
\tau_S : H_1(T_S) \to L_\Z/H_\Z.
$$
A consequence of Johnson's Theorem (Thm.~\ref{thm:johnson}) is:

\begin{corollary}
If $g\ge 3$, then $\tau_S$ is surjective with finite kernel of exponent $2$.
\qed
\end{corollary}

\subsection{Cohomology computations}

Arguments similar to those in Section~\ref{sec:monod_comps} can be used to
prove the following result.

\begin{proposition}
If $g\ge 3$, then
\begin{enumerate}
\item If $V$ is an irreducible rational representation of $\Sp(H_\Q)$ that
is not isomorphic to $L_\Q/H_\Q$, then $H^1(\G_S,V)$ vanishes.
\item The homomorphism
$$
H^1(\G_S,L_\Z/H_\Z) \to H^1(\G_{S,x},L_\Z/H_\Z)
$$
induced by the canonical surjection $\G_{S,x} \to \G_S$ is an isomorphism.
\item The group $H^1(\G_S,L_\Z/H_\Z)$ is freely generated over $\Z$ by a class
whose image under the restriction mapping
$$
H^1(\G_S,L_\Z/H_\Z) \to \Hom_{\Sp(H_A)}(H_1(T_S),L_\Z/H_\Z)
$$
is $2\tau_S$.
\end{enumerate}
\end{proposition}

Note that the \'etale local system $\L_\Zl/\H_\Zl$ over $\M_g^1$ is actually
the pullback of a local system over $\M_g$ for which we use the same notation.
Specifically, this is the local system over $\M_g$ corresponding to the third
fundamental representation of $\GSp(H_\Ql)$.

\begin{corollary}
If $g\ge 3$, then the mapping
$$
\Het^1(\M_g,\L_\Zl/\H_\Zl) \to \Het^1(\M_g^1,\L_\Zl/\H_\Zl)
$$
induced by the projection $\M_g^1 \to \M_g$ is an isomorphism. If $\V$ is the
\'etale local system over $\M_g$ corresponding to an irreducible 
representation of $\GSp(H_\Ql)$ that is not isomorphic to $\L_\Ql/\H_\Ql$ or
some $\Ql(n)$, then $\Het^1(\M_g,\V)$ vanishes. In addition the structure
morphism $\M_g \to \Spec K$ induces an isomorphism
$$
H^1(G_K,\Ql(n)) \cong \Het^1(\M_g,\Ql(n)). \qed
$$
\end{corollary}

We shall denote the class in $\Het^1(\M_g,\L_\Zl/\H_\Zl)$ that corresponds to
$\nu(\cC_g^1)\in \Het^1(\M_g^1,\L_\Zl/\H_\Zl)$ by $\nu(\cC_g)$. Note that
$\Het^1(\M_g,\L_\Zl/\H_\Zl)$ is freely generated by it. This will allow us to
extend, in the next section, the definition of $\nu(C)$ given in
Section~\ref{sec:nu} to curves $C\to B$ without a $B$-rational point.

\subsection{The invariant $\nu(C)$}
\label{sec:nu}

The computations of the previous section can be used to define the invariant
$\nu(C) \in \Het^1(B,\L_\Zl/\H_\Zl)$ even when $C$ does not have any rational
points, which proves Proposition~\ref{prop:equal}.

Suppose that $C\to B$ is a family of smooth projective curves of genus $\ge 3$
over $K$. It is classified by a morphism $[C] : B \to \M_g$.
Define
$$
\nu(C) \in \Het^1(B,\L_\Zl/\H_\Zl)
$$
to be the pullback $[C]^\ast \nu(\cC_g)$ of the universal class. It follows
from the results in the previous section that
$$
\nu(C) = \nutilde(C,\xi)
$$
when $C$ has a $B$-rational point $\xi$.

Theorem~\ref{th:main-out} can now be proved using arguments similar to those in
Section~\ref{sec:density}.

\subsection{An un-pointed Harris-Pulte Theorem}
\label{sec:harris-pulte_unptd}

There is a version of the $\l$-adic Harris-Pulte Theorem for curves without
any rational points.

Let $n(C)\in H^1(G_K,L_\Zl/H_\Zl)$ be the pullback of the universal class
$$
\nhat^\prol \in H^1(\Out \pi_1^\prol(C^\an), L_\Zl/H_\Zl)
$$
along the homomorphism $G_K \to \Out \pi_1^\prol(C^\an)$.

\begin{theorem}
\label{th:unptd-Harris-Pulte}
If the genus of $C$ is $\ge 3$, then the classes $n(C)$ and $\nu(C)$ are
equal in $H^1(G_K, L_\Zl/H_\Zl)$. \qed
\end{theorem}

\subsection{The $\l$-adic formulation in the un-pointed case}

As in previous sections, $C$ is a smooth projective curve of genus $\ge 3$
defined over the subfield $K$ of $\C$, $\Kbar$ is its algebraic closure in
$\C$, $x$ is the geometric point $\Spec \Kbar$.

We have the truncations
$$
\rho_C^{\prol,m} : G_K \to \Out \big[\pi_1^\prol(C^\an)/L^{m+1}\big]
$$
and
$$
\rhohat_C^{\prol,m} : \pi_1^\alg(\M_g,[C\otimes\C]) \to
\Out \big[\pi_1^\prol(C^\an)/L^{m+1}\big]
$$
of the monodromy representations $\rho_C^\prol$ and $\theta_C^\prol$
defined in Section~\ref{sec:monod}.


\renewcommand{\theenumi}{\alph{enumi}}

\begin{theorem}
\label{th:main-out-4}
If the image of the $\l$-adic cyclotomic character $\chi_\l : G_K \to
\Zl^\times$ is infinite, then the following four conditions are equivalent:
\begin{enumerate}
\item $\im\rho_C^{\prol,m}$ is an open subgroup of $\im\rhohat_C^{\prol,m}$ for
all $m\ge 1$.
\item $\im\rho_C^{\prol,2}$ is an open subgroup of $\im\rhohat_C^{\prol,2}$.
\item The image of $\rho^{\prol,1}_C : G_K \to \GSp(H_\Zl)$ is open, and the
class $\nu(C)$ in $H^1(G_K,L_\Zl/H_\Zl)$ has infinite order.
\item Any of the 3 equivalent conditions (i), (ii), (iii) of
Theorem~\ref{th:main-out}. \qed
\end{enumerate}
\end{theorem}

\bigskip

\noindent{\it Acknowledgments.} We would like to thank A.\ Tamagawa, S.\
Mochizuki, and K.\ Fujiwara for helpful conversations related to this work.
We are also grateful to N.~Kawazumi
and S.~Morita for correspondence related to the Magnus homomorphism.

\end{document}